%

\input amstex
\magnification=1200
\loadmsam
\loadmsbm
\loadeufm
\loadeusm
\UseAMSsymbols

\hsize=6.9truein
\hoffset=-0.11truein
\vsize=8.9truein
\voffset=-0.2truein

\def\leftitem#1{\item{\hbox to\parindent{\enspace#1\hfill}}}

\def\boxit#1#2{\hbox{\vrule
	\vtop{%
	\vbox{\hrule\kern#1%
	\hbox{\kern#1#2\kern#1}}%
	\kern#1\hrule}%
	\vrule}}

\def\leaderfill{\leaders\hbox to 1em{\hss.\hss}\hfill}

\parskip=\medskipamount
\document

\input epsf



\centerline{\bf Characters of SL(2) Representations of Groups}

\medskip
\centerline{\bf Feng Luo}

{\bf \S1. Introduction}

1.1. Given a field $K$ and a representation $\rho$ of a group to $SL(2,K)$,
the \it character \rm of the representation sends a group element $g$ to
the trace of the matrix $\rho (g)$. 
One of the result of the paper  is the following,

{\bf Theorem.} \it
Suppose $K$ is a field  
so that each quadratic equation with  coefficients in $K$ has a root in
$K$.
Then a $K$-valued function defined on a
group  is the character of an $SL(2,K)$ representation of the 
group if and only if 
its restriction to each 2-generator subgroup is an $SL(2,K)$ character.
\rm

The $SL(2,K)$ characters on 2-generator groups are well understood since
the work of  Fricke-Klein [FK] and Vogt [Vo]. They are governed by the
trace identity: $tr(AB) + tr(A^{-1} B) = tr(A) tr(B)$ for $SL(2,K)$
matrices $A, B$. In [Hel], Helling gave an elegant axiomatic  approach
to characters based on the above trace identity. Following Helling,
a $K$-valued  function $f$ defined on a group
$G$ is called a \it $K$-trace function \rm 
if (1) $f(xy) + f(x^{-1}y) = f(x) f(y)$ for all $x,y$ in $G$ and
(2) $f(id) =2$ where $id$ is the identity element. 
Evidently an $SL(2,K)$ character is a $K$-trace
function. Using the work of
[FK], [Hel] and [Vo] that $K$-trace functions are characters on
2-generator groups, one deduces  the following result equivalent
to the above theorem. 

{\bf Corollary.} \it 
Suppose $K$ is a field
so that each quadratic equation with  coefficients in $K$ has a root in
$K$. Then a $K$-valued function $f$ defined on
a group is the character of an $SL(2,K)$ representation  of the group if
and only if  $f(xy) + f(x^{-1}y) = f(x) f(y)$ for all
$x,y$ in the group and $f(id) =2$. \rm

This generalizes a result of Helling [Hel] who proved that
$\bold R$-trace functions are $SL(2,\bold R)$ characters  under some
additional assumptions.


1.2. The main result of the paper which implies theorem 1.1 
gives a characterization of $SL(2,K)$ 
characters defined on the fundamental groups of surfaces using subsurface
groups.  
In contrast to theorem 1.1 which  uses
the hierarchy of subgroups indexed by the number of generators to
describe the characters,  there exists a natural hierarchy
of surfaces under inclusion indexed by the level.
Recall that the \it
level \rm of 
a compact surface of negative Euler number is the minimal number
of disjoint simple loops decomposing
the surface into 3-holed spheres.  It is also the complex dimension
of the Teichm\"uller space of complex structures on the interior of
the surface with punctured ends. For instance, the 3-holed sphere
has level-0 and  the 4-holed sphere and  the 1-holed torus have level-1.
This  hierarchy of surfaces is prominent in Grothendieck's manuscript
[Gr] and conformal field theory [MS]. In particular, 
 Grothendieck conjectured that  the ``tower 
of Teichm\"uller spaces" can be reconstructed  from the  Teichm\"uller
spaces of level-1 surfaces subject to the relations supported in 
level-2 surfaces. Motivated by this Grothendieck's reconstruction 
principle, one asks if a character can be reconstructed from its 
restriction to the fundamental group of each level-1 subsurface. 
The main result of the paper  gives a complete answer to this question.

To give a precise solution to the above question, we shall first
note that the fundamental group is not vital with respect
to the hierarchy. Indeed, if a given element in the  fundamental 
group is ``complicated" in the sense that it has no representative 
in any level-1 subsurface, then the condition becomes null about the class. 
Thus, we focus our attention to those ``simple" elements in the fundamental
group.  This motivates the introduction of
the set $S(\Sigma)$ of free
homotopy classes of unoriented homotopically non-trivial simple
loops on a surface $\Sigma$. The space $S(\Sigma)$ was 
introduced by Max Dehn [De] in  his study of the mapping class groups and
was independently introduced by Thurston [Th] in his work
on surface theory.
If $f$ is an $SL(2,K)$ character defined
on the fundamental group of a surface $\Sigma$, 
then $f$ induces a $K$-valued
function on $S(\Sigma)$ which we still call an $SL(2,K)$ character.
A natural property of a character $f$ on $S(\Sigma)$ is that its
restriction to each subsurface is again a character. To be more
precisely, if $\Sigma'$ is an essential subsurface (i.e.,
the inclusion map induces a monomorphism between fundamental groups),
then the \it  restriction map \rm  $f \circ i_*$ is again a character on
$S(\Sigma')$ where $i_*:$ $ S(\Sigma') \to S(\Sigma)$ is induced by
the inclusion. 

We call a $K$-valued function $f$ defined
on the set $S(\Sigma)$ of homotopy classes of simple loops a
\it trace function \rm if the restriction of the function
to each $S(\Sigma')$ is a character for each level-1 essential
subsurfaces $\Sigma'$.  Grothendieck's principle predicts that a
trace function is a character. The main result of the 
paper shows that this holds except for finitely many exceptional 
trace functions defined on genus zero surfaces when the characteristic
of the field $K$ is not 2. 
All exceptional trace functions are derived from a single
one defined on the 5-holed sphere which we describe as follows.
Let the characteristic of the field $K$ be not 2 and
$b_1, ..., b_5$ be the boundary components of the 5-holed sphere
$\Sigma_{0,5}$. Define $f_0: S(\Sigma_{0,5}) \to K$ by sending  each
$b_i$ to 2 and all other elements to $-2$. One checks easily 
(see \S5.4) that $f_0$ is a trace function which is not
the character of any representations. There are
 sixteen exceptional trace functions $f$ on the 5-holed sphere
all derived from $f_0$. Namely,
an \it exceptional trace function \rm $f: S(\Sigma_{0,5}) \to   K$ satisfies the
following (1) $f(S(\Sigma_{0,5})) = \{2, -2\}$, (2) $\Pi_{i=1}^5
f(b_i) = 32$, and (3) if $\alpha$ is a non-boundary parallel class
so that $\alpha, b_i$, and $b_j$ bound a 3-holed sphere, then $f(\alpha)
= -\frac{1}{2} f(b_i) f(b_j)$. 

The main result of the paper is the following.

{\bf Theorem.} \it 
Suppose $K$ is a field
so that each quadratic equation with  coefficients in $K$ has a root in
$K$. Let $f$ be a $K$-valued trace function defined on the set
 $S(\Sigma)$ of homotopy
classes of essential simple loops in a compact orientable surface $\Sigma$.

(1) If the characteristic of the field $K$ is 2, then the 
trace function $f$ is the character of an $SL(2,K)$ representation.

(2) If the characteristic of the  field $K$ is not 2, then

(2.1) either $f$ is the character of an $SL(2,K)$ representation, or

(2.2) the genus of the surface $\Sigma$ is zero, $f$ takes
only values $\{2, -2\}$, and there is a level-2 
subsurface so that the restriction of $f$ to the subsurface
is one of the  sixteen  exceptional trace functions.

(2.3) There exist
exceptional trace functions on each genus zero surface of level
at least $2$. The the number of
exceptional trace functions  on a fixed surface is finite.

\rm

Note that surfaces in the theorem are connected and could be compact
or non-compact of infinite type. 

Theorem 1.2 does not cover the case when the surface $\Sigma$
has level at most 1.
The characterization of $SL(2,K)$ characters on the set of
simple loops in level-1 surfaces (propositions 3.4 and 3.5) 
is well known by the work of [FK],  [Go], [Hel], [Ho], [Ma], [Vo] and others. 
It is based on the following lemma  (lemma 2.3) well known
to the experts in the field.
Namely, given six elements $x_1, x_2, x_3,
x_{12}, x_{23}, x_{31}$ in $K$, there exist three $SL(2,K)$ matrices
$A_1, A_2,$ an $  A_3$ so that $tr(A_i) = x_i$ and $tr(A_i A_j) = x_{ij}$.

1.3.  Given a group $G$ and a field $K$, the set of all $SL(2,K)$ characters
on $G$ is called the \it character variety \rm
of the group.
Theorem 1.1 gives an explicit algebraic description of the
character variety of the group for those field $K$ satisfying the
condition in the theorem. If the group $G$ is finitely generated,
then a well known result (proposition 2.2)
shows that there
exists a finite subset $F \subset G$ so that each $SL(2,K)$ character
on the group is algebraically determined by its restriction to the
finite set $F$.  
As a  consequence of these and  the Hilbert basis theorem, one
obtains the following corollary which slightly generalizes  
a result of Culler-Shalen  [CS] who proved it for algebraically closed
field $K$.

{\bf Corollary}(Culler-Shalen).  \it
Suppose $G$  is a finitely generated group and $K$ is a field 
so that each quadratic equation with  coefficients in $K$ has a root in
$K$. Then the  set of all $SL(2,K)$ characters on the group
forms an affine algebraic variety defined over $K$.
Furthermore, the defining equations are integer coefficient
polynomials. \rm

In [GM], Gonz\`alez-Acu\~na and Montesinos-Amilibia
gave a constructive proof of Culler-Shalen's result. Their
proof also shows the above corollary in the case the characteristic of
$K$ is not $2$.

The assumption on the quadratic closeness of the field $K$
can be replaced by extension fields. Namely, suppose $K$ is any field and 
$f$ is a $K$-trace function defined on a group
generated by $n$ elements. Then there exists an  extension field
$F$ of $K$ obtained from $K$ by at most $n$ quadratic extensions
and a representation
of the group to $SL(2, F)$ whose character is the given $K$-trace function.


1.4. There exists an interesting analogy between the hierarchy
of finitely  generated  groups indexed by the number of generators
and the hierarchy of surfaces indexed by the level. It seems
that 
the role of level-1 surfaces is similar to that of
2-generator groups. For instance, Jorgensen [Jo] proved
that a non-elementary subgroup of $SL(2,\bold C)$ is discrete if and 
only if each 2-generator subgroup is discrete. A consequence of [Lu1] shows that
a faithful representation of a surface group to $SL(2,\bold R)$ is
discrete if and only if the restriction of the representation to 
each level-1 subsurface group is discrete and uniformizes the subsurface.
Theorems 1.1 and 1.2 provide another comparison. 
Here is a third pair. Recall that a subgroup in $SL(2,K)$ is \it
reducible \rm if it leaves a 1-dimensional 
linear subspace in $K^2$ invariant. It is known [CS] 
that a subgroup in $SL(2,K)$ is reducible if and
only if each 2-generator subgroup is reducible (see \S2.5). The analogous
result is the following. 

{\bf Theorem.}  \it An $SL(2,K)$ representation of a
surface group is reducible if and only if
its restriction
to each level-1 subsurface group is reducible.
\rm

In fact,  in the statement of the theorem,   3-holed
sphere and 1-holed torus subgroups suffice. 
However, there exists an irreducible representation of a surface group to
$SL(2, \bold K)$ so that the restriction
to each level-0 subsurface group is reducible. Such irreducible
representations  occur  rarely (only on genus 1 surfaces)
and  are classified in \S6 and \S7.

Finally, the analogous result to the well known 
 proposition 2.2 is
the following  (see \S3.9) that
\it there exists a finite set of homotopy classes of
 simple loops on each compact orientable  surface
so that the characters of $SL(2, K)$ representations are algebraically
determined by their restrictions to the finite set. \rm

1.5.
Since each compact 3-manifold has a Heegaard splitting, a 3-manifold group
is the quotient of a surface group by a subgroup of the form $N_1 N_2$
where each $N_i$ is  normally generated by disjoint simple loops. This shows
that simple loops are characteristic for 3-manifold groups (among all
finitely presented group). By singling out the special  feature of
simple loops in theorem 1.2, it is  hoped that it will have some applications
to 3-manifold groups.  In particular, we are motivated by the following question.
Given a Haken 3-manifold $M$, does there exist an irreducible $SL(2, K)$
representation of the fundamental group of the 3-manifold for some finite
field $K$? See [Hem] for related topics.

1.6. As mentioned before, theorem 1.2  may be interpreted as 
establishing  Grothendieck's
reconstruction principle  for $SL(2)$ character varieties.  Broadly
speaking, the principle says that to  study the isotopy
class of a structure on a surface, one should consider the restriction
of the structure to  the isotopy classes of all level-1 subsurfaces
and reconstruct the original isotopy class of the  structure from  the
restrictions. Furthermore, level-2 subsurfaces should serve as the
``relators" in the reconstruction process (see [Gr] and [Lu4]).
This reconstruction principle is shown to be valid
for hyperbolic metrics and measured laminations in [Lu1] and [Lu2].
The proof of theorem 1.2 is similar to the proof  of [Lu1]. Namely,
first we prove the result for level-1 and level-2 surfaces and
then we prove the result for all surfaces using a general gluing
lemma.  The main
difficulty in proving theorem 1.2 is caused by the existence of
irreducible representations whose restrictions to some 2-generator subgroups
are reducible.
Similarly, the main difficulty in establishing theorem 1.1 is
in the case of free group on 4 generators. 

1.7. The study of the algebra of characters of $SL(2,K)$ representation
was started by Vogt and Fricke-Klein and is developed by many authors 
[BG], [BH], [CS], [Hel], [Ho], [Ke],  [LM], [Ma], [Pr], [Sa] and others. 
It seems that there is a close relation between
what we did here and those algebraic approach to the ring of
$SL(2)$ characters on a group.

1.8. The organization of the paper is as follows.
In \S2, we recall the basic facts on traces of $SL(2,K)$ matrices and
group representations.
In \S3, we recall the basic facts on simple loops on surfaces and the
modular structure.  We then use the modular structure to describe
$SL(2,K)$ characters on level-1 surfaces. A multiplication of the
simple loops on surfaces will also be discussed. In \S4 and \S5, we prove
the main result for the genus zero surfaces  by making
extensive use of the modular structure. In
\S6, we  prove theorem 1.2 for the 2-holed torus.
Theorem 1.2
for all surfaces is proved in \S7. In \S8, we prove theorem 1.1.
In the final section \S9, 
we discuss some questions arising from the consideration of $SL(2)$ characters.

1.9. \it Acknowledgment. \rm I would like to thank F. Bonahon and
 X.-S. Lin for many
discussions. This work is supported in part by the NSF.

\S2. {\bf Preliminaries on SL(2,K) Matrices}

In this section, we shall introduce notations and  recall basic trace 
identities.

2.1. We shall use the following notations and terminologies. 
A representation $\rho$ of a group $G$ to $SL(2,K)$ is called 
\it reducible \rm if there exists a 1-dimensional linear subspace 
in $K^2$ invariant under
the linear action of $G$. Otherwise the  representation  is called \it
irreducible\rm. 
Two
representations $\rho_1$ and $\rho_2$ of a group $G$ to
 $SL(2,K)$ are \it conjugate \rm if there exists a matrix $X$ in  $SL(2,K)$ so that
for all $g \in G$, $X\rho_1(g) X^{-1} = \rho_2(g)$. Evidently, conjugate
representations have the same characters.
 A reducible representation is called \it diagonalizable
\rm if it is $SL(2,K)$ conjugate to a representation whose image
lies in the  set of diagonal matrices. 
The character of a reducible (resp. irreducible)
representation is also called reducible (resp. irreducible). 
A subgroup
of $SL(2,K)$ is called  reducible if the inclusion map is reducible.

2.2.  The following trace identities will be used  frequently.  They
are derived from the first identity (a). The earliest source of
these identities seems to be [Vo]. For instance, the less commonly used
identity (d)  is on page S11 in [Vo]. See [FK], [Go],
[Vo] and others for a proof.

{\bf Lemma.} \it  Let $R$ be a commutative ring with identity.
Suppose $A, B,  A_i, B_i$ are $SL(2,R)$ matrices  where $i=1,2,3.$
Then the following identities hold.

(a) $tr(AB) + tr(A^{-1}B) = tr(A) tr(B).$

(b)  $tr^2(A) + tr^2(B) + tr^2(AB) - tr(A) tr(B) tr(AB) = tr([A,B]) +2.$

(c) Let $A_4=  A_1$, $A_5 = A_2$, $P = \sum_{i=1}^3 tr(A_i) tr(A_{i+1}
A_{i+2}) - tr(A_1) tr(A_2) tr(A_3)$ and 
$Q = \sum_{i=1}^3 (tr^2(A_i) + tr^2(A_i A_{i+1}) - tr(A_i) tr(A_{i+1})
tr(A_{i}A_{i+1})) + tr(A_1A_2) tr(A_2A_3) tr(A_3A_1) -4$. Then the two
roots of the quadratic equation $x^2 -Px + Q = 0$ are
$tr(A_1A_2A_3)$ and $tr(A_1^{-1} A_2^{-1} A_3^{-1})$.

(d) $tr(A_1A_3) + tr(A_1 A_2 A_3 A_2^{-1}) = - tr(A_1A_2) tr(A_2A_3)
+ tr(A_1) tr(A_3) + tr(A_2) tr(A_1A_2A_3).$ \rm


In \S3, these equations will be interpreted using the
modular configuration  $(\hat \bold Q, PSL(2, \bold Z))$.


As a consequence of the lemma, one has the following useful proposition.
See [CS], [FK], [Ho] and [Vo] for a proof.

{\bf Proposition} ([CS], [FK], [Ho] and [Vo]). \it
Given  a commutative ring $R$ with identity, 
the trace of a word $w(A_1, ..., A_n)$ 
in the $SL(2, R)$ matrices $A_1, ..., A_n$
is a polynomial with integer coefficients in the
traces of  $A_{i_1} ...$$A_{i_k}$ where $1 \leq i_1 < ...< i_k \leq n$
and $k \leq n$. 
 \rm


 \it Remark. \rm
If the commutative ring $R$ is a field $K$ of characteristic not
equal to 2, then a  trace identity in 
[Vo] (page S14, line 19) shows that the
trace $tr(A_1A_2A_3A_4)$ can be expressed in terms of the traces of
 $A_i,  A_iA_j, $ and $A_iA_jA_k$, $1 \leq i,j,k \leq 4$. 
Thus in this case, one can strengthen
the proposition to $tr(A_{i_1} ... A_{i_k})$ where $1 \leq i_1 < ...< i_k
\leq n$ for $k \leq 3$. This triple-trace theorem has been rediscovered
independently by many mathematicians. See [Bu], [BG] and others.
 
In [Hel], Helling proved that all trace identities in Lemma 2.2 still
hold for $R$-trace  functions. Since the proof of the above proposition
uses only identity (a) in lemma 2.2 and $tr(id) = 2$, Helling proved 
the following corresponding result for $R$-trace functions.

{\bf Corollary}([Hel]). \it Suppose $G$ is a group generated by
$n$ elements $\{x_1, ..., x_n\}$.  Then for each element
$w \in G$, there exists an integer coefficient polynomial
$P_w$ in variables $t_{i_1 ... i_k}$, where   $1 \leq i_1 < ...< i_k
\leq n$ for $k \leq n$ so that for all $R$-trace functions $f$ on $G$,
$f(w) = P_w( f(x_1), ..., f(x_{i_1} ... x_{i_k}), ..., f(x_1... x_n))$.
\rm

\midspace{0.1cm}
\centerline{\epsfbox{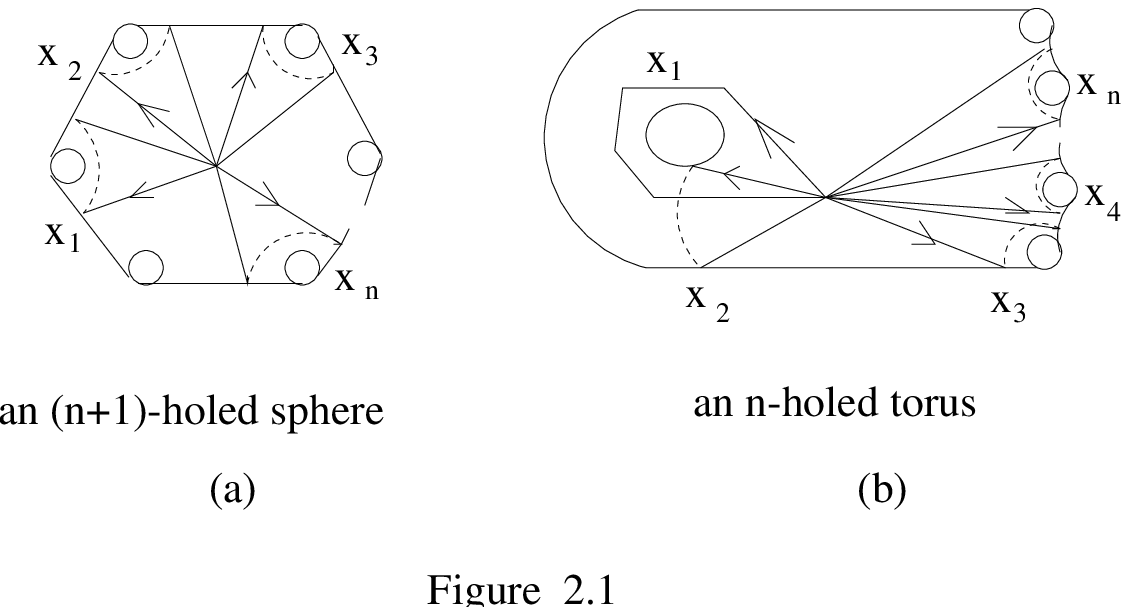}}

In particular $R$-trace functions on the free group $F_n$ of
$n$ generators $<x_1, ..., x_n>$ are determined by their
restrictions to the set $\{x_{i_1} ... x_{i_k} |$ 
$1 \leq i_1 < ...< i_k \leq n$ and $k \leq n$\}.
Now identify the free group $F_n$ with the fundamental group
$\pi_1(\Sigma_{0, n+1})$ of the $(n+1)$-holed sphere or
$\pi_1(\Sigma_{1, n-1})$ of the $(n-1)$-holed torus. Then
we may choose the set of $n$ generators $\{x_1, ..., x_n\}$ 
so that each element $x_{i_1} ... x_{i_k}$
is represented by a simple loop in the surface (see fig. 2.1.). This
shows that $R$-trace functions on the  fundamental groups of
these surfaces are determined by their restrictions to the
classes of simple loops. 

For the low-rank free groups
 $F_2, F_3$ and $F_4$, the $(2^{n}-1)$-element set 
 $\{x_{i_1} ... x_{i_k} |$ $1 \leq i_1 < ...< i_k \leq n$\}
 are closely related to
the so called modular relation and pentagon relations on
the surfaces of genus zero. 
Indeed, take $F_2 = \pi_1(\Sigma_{0,3})$, then the 3-element set
$\{x_1, x_2, x_1 x_2\}$ is represented by the three boundary
components.  Take $F_3 =\pi_1(\Sigma_{0,4})$. Then the 7-element
set $\{ x_1, x_2, x_3, x_1 x_2, x_2 x_3, x_1 x_3, x_1 x_2 x_3\}$
 is represented by the four boundary components and three simple loops
pairwise intersecting at two points (see fig. 2.2). 
Take $F_4 = \pi_1(\Sigma_{0,5})$. Then the 15-element set
$\{x_1, x_2, x_3, x_4, x_1x_2, x_1x_3, x_1x_4, x_2 x_3,$
$ x_2x_4,$  $
 x_3x_4, x_1x_2x_3,$$ x_1x_3 x_4,$$ x_2x_3x_4, x_1x_2x_4,$
$ x_1x_2x_3x_4\}$ is 
represented by the five boundary components and 10 more simple loops
closely related to the pentagon relation (see \S3.2 and \S4
 for more discussions).

\midspace{0.1cm}
\centerline{\epsfbox{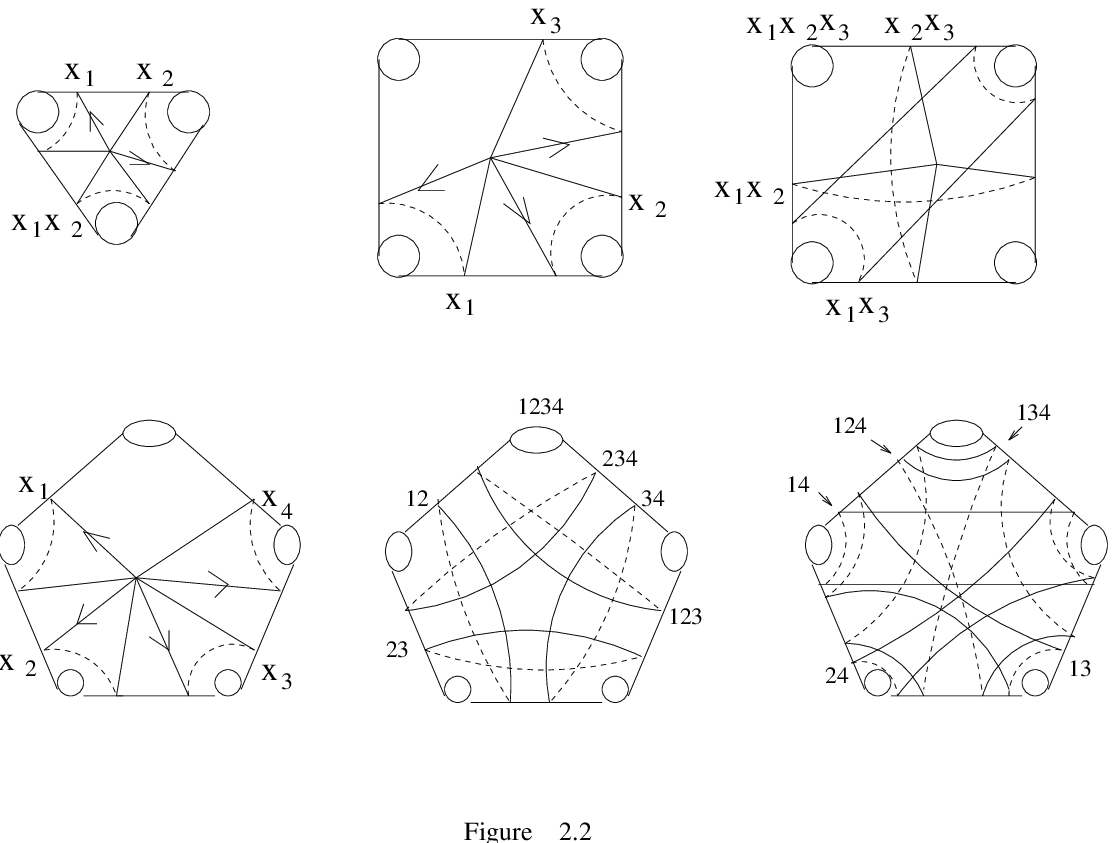}}

2.3.
In the rest of the paper, we will always assume that the field
$K$ is \it quadratically closed \rm in the sense that each quadratic
equation $x^2 + ax + b = 0$, $a,b \in K$ has roots in K.
Under this assumption, each $SL(2,K)$ matrix has
eigenvalues in $K$.
Furthermore, two $SL(2,K)$ matrices are $SL(2,K)$ conjugate if
and only if they have the same trace.

The following lemma is known to experts in the field. Especially,
in the case the field $K$ is algebraically closed, it follows
from the combination of the work of Culler-Shalen [CS] and 
Horowitz  ([Ho], theorem 4.3).
Since we have not
seen a written proof of the version stated below, a  proof is given 
in the appendix for completeness.

{\bf Lemma.}  \it Suppose $K$ is a quadratically closed field.
Given six elements $t_1, t_2, t_3, t_{12}, t_{23}$
and $t_{31}$ in $K$, there exist three $SL(2,K)$ matrices $A_1, A_2$ and
$A_3$ so that $tr(A_i) = t_i$ and $tr(A_i A_j)  = t_{ij}$. \rm

Combining the lemma with proposition 2.2 and  lemma 2.2(c),
one sees that the set of all characters on
the free group in three generators $<x_1, x_2, x_3>$  is
the hypersurface $\{(t_1, t_2, t_3, t_{12},$ \newline $t_{23},$$ t_{31},$
$ t_{123})
\in K^7 |$ the equation $(*)$ holds\},
$$t_{123}^2 +(\prod_{i=1}^3 t_i - \sum t_i t_{jk}) t_{123}
+ \sum_{i=1}^3 t^2_i + \sum t^2_{ij} + \prod t_{ij} - \sum t_i t_j t_{ij} -4 =0.
\tag *
$$
In the equation $(*)$, $t_{i_1...i_k} = tr (\rho(x_{i_1} ... x_{i_k}))$
and the indices $i,j,k$ are  pairwise distinct. This fact was well known
after the work of Culler-Shalen [CS] and Horowitz [Ho].
See \S3.5 for more details. 

2.4. It is easy to see that  $(BA_1B^{-1}, BA_2B^{-1}, BA_3B^{-1})$
and $(A_1^{-1}, A_2^{-1}, A_3^{-1})$ are other solutions in lemma 2.3.
In fact, these are the set of all solutions to the the
equation $tr(X_i) = t_i, tr(X_iX_j ) = t_{ij}$ if and only if the group
generated by $<A_1, A_2, A_3>$ is irreducible. To derive this, let us
recall  the 
following lemma proved by Culler and Shalen ([CS] lemma 1.5.2).

{\bf Lemma} (Culler-Shalen). \it Suppose the field $K$ is
quadratically closed.
If $\rho_1$ and $\rho_2$ are two representations of a group to $SL(2,K)$ 
so that $\rho_1$ is irreducible, then $\rho_1$ is conjugate
to $\rho_2$ if and only  if they have the same character functions.\rm

As a consequence of Culler-Shalen's lemma, lemma 2.2 and proposition 2.2,
we see that if the
group generated by $\{A_1, A_2, A_3\}$ in lemma 2.3 is irreducible,
then the solution
$(A_1, A_2, A_3)$ in lemma 2.3 is unique up to conjugation and inverse.
Evidently, if the group $<A_1, A_2, A_3>$ in lemma 2.3 is reducible,
then the solution is not unique in the above sense.

2.5.  Due to the importance of irreducible representations, we need 
an irreducibility criterion. The following is well known. See [CS]
lemma 1.5.5, or [Ma] for instance. 

{\bf Lemma.} \it Suppose the field $K$ is quadratically closed.
 The group $<A, B>$ generated by two elements $A,B$ in $SL(2,K)$
is reducible if and  only  if $tr([A,B]) = 2$. \rm

By lemma 2.2(b), the condition in the above lemma is the same as
$tr^2(A) + tr^2(B) + tr^2(AB) - tr(A) tr(B) tr(AB) -4 =0$. Since this
expression will occur frequently, following [Ma], let us denote 
$tr^2(A) + tr^2(B) + tr^2(AB) - tr(A) tr(B) tr(AB) -4$ by $\Delta (A,B)$.

As a consequence of the lemma, one  has the following criterion 
for reducibility. A 
slightly different criterion can be found in [CS],
[GM] (proposition 4.4) and [Ma].

{\bf Corollary.} \it  Suppose the field $K$ is quadratically closed.

(a).  The group $<A_1, A_2, A_3>$ generated by three elements $A_1, A_2$
and $A_3$ in $SL(2,K)$ is reducible if and only if $\Delta (A_i, A_j) =0$ and
$\Delta (A_1, A_1A_2A_3) = 0$ where $(i,j) =(1,2), (2,3), (3,1)$.

(b). The group $<A_1, ..., A_n>$ in $SL(2,K)$ generated  by  $n$ elements
is reducible if and only if  each 3-generator subgroup $<A_i, A_j, A_k>$
is reducible, i.e., for all possible choice of indices $i,j,k$,
$\Delta (A_i, A_j) = 0$ and $\Delta (A_i, A_i A_j A_k) = 0$.

(c) (Culler-Shalen). A subgroup of $SL(2,K)$ is reducible if and only if 
each 2-generator subgroup is reducible. 

\rm

\it Proof. \rm  To see part (a), we may assume that none of $A_i$ is
$\pm id$ since otherwise it reduces to  lemma  2.5.  Thus each $A_i$ has
at most two eigenspaces. By lemma 2.5, each pair $(A_i, A_j)$ has 
a common eigenspace for $i \neq j$. Now if one of $A_i$ has exactly
one eigenspace, then all $A_1, A_2, $ and $A_3$ share this unique
eigenspace. Thus the group is reducible. 
If otherwise, each $A_i$ has two distinct eigenspaces $L_j$
and $L_k$, $i \neq j \neq k \neq i$.  Now suppose that  the group 
$<A_1, A_2, A_3>$ is irreducible.  
Then all three eigenspaces are pairwise distinct,
i.e., $L_1 \neq L_2 \neq L_3 \neq L_1$. But by assumption, $A_1$ and
$A_1 A_2 A_3$ have a common eigenspace $L$. Due to $A_1(L) = L$,
thus  $L$ must be either $L_2$ or $L_3$, say,  $L = L_2$. Then 
$A_1A_2A_3(L_2) = L_2$ implies that $A_2(L_2) = L_2$, i.e., $L_2$ is either $L_1$ or $L_3$
which contradicts the assumption.

Parts (b) and (c) follow from part (a) easily. To see (b), we first drop
all generators $A_i$ which are $\pm id$. Thus, we may assume that
each $A_i$ has at most two eigenspaces.  By the assumption, any three
elements $A_i, A_j,$ and $A_k$ have a common eigenspace. The goal is to
show that all elements $A_i$ have a common eigenspace. To this end, we form
a graph whose vertices are eigenspaces of $A_i$'s. To each element $A_i$,
we draw an edge ending at the eigenspaces of $A_i$ (the edge becomes a
loop if $A_i$ has only one eigenspace).  Now by the assumption, any three
edges of the graph has a common vertex. Thus all edges of the graph share
a  vertex.

To see part (c), take a subgroup $G$ of $SL(2,K)$  which has the
property that each 2-generator subgroup is reducible. By parts (a) and
(b), we see
that each finitely generated subgroup of $G$ is reducible. Now by the
same graph theoretical argument as in part (b), we see that the group G
is reducible.  $\square$


2.6. Given a reducible representation $\rho$ of a group to
$SL(2, K)$, we may assume after conjugate $\rho$ by an $SL(2,K)$
matrix that the image of $\rho$ is in the set of upper triangular matrices.
Let $\rho'$ be a new representation so that $\rho'(g)$ is the
diagonal matrix whose diagonal entries are that of $\rho(g)$.
Then the diagonalizable representation $\rho'$ has the same
character as that of $\rho$. Evidently the diagonal representation
is unique up to conjugation.  We call $\rho'$ the \it diagonalization
\rm of the reducible representation $\rho$. 

{\bf Lemma.} \it Suppose $K$ is a quadratically closed field. Let
$\rho_1$ and $\rho_2$ be two non-diagonalizable
reducible
representations of the free group $<x,y>$ on two generators to $SL(2,K)$. 
If $\rho_1$ and $\rho_2$ have the same character and $tr(\rho_i(x)) 
\neq \pm 2$, then
they are conjugate. \rm

Indeed, under the assumption, we can  conjugate the 
pair ($\rho_i(x), \rho_i(y)$)
to the pair of  matrices  ($\left(\matrix \lambda & 0\\ 0 & \lambda^{-1}
 \endmatrix \right)$,  $\left(\matrix \mu & 1\\ 0 & \mu^{-1}
\endmatrix \right)$) where $\lambda \neq \pm 1$.

{\bf \S3. Simple Loops on Surfaces and the  Modular Configuration}

We shall recall some basic facts on the set 
of isotopy classes
of simple loops on surfaces
and  express the results in \S2 in terms of a $(\bold QP^1, PSL(2, \bold Z))$
modular structure on the set $S(\Sigma)$.
We also establish several irreducible conditions in terms of
the modular structure. 

The field $K$ is always assumed to be quadratically closed.

3.1.  The following notations and terminologies  will be used. Let 
$\Sigma$ = $\Sigma_{g,r}$ be a compact orientable surface of genus  $g$ 
with $r$ boundary components.  The \it level \rm 
of the  surface $\Sigma_{g,r}$  is defined to be $3g+r-3$ which
is the minimal number of disjoint simple
loops decomposing  the surface into 3-holed spheres. 
Recall that $ S(\Sigma)$ is  the set of isotopy (homotopy) 
classes of unoriented homotopically non-trivial simple loops on 
$\Sigma$. Let $ S'(\Sigma)$ be the subset of $ S(\Sigma)$ consisting
of non-boundary parallel isotopy classes. 
The fundamental group of the surface is denoted by
$\pi_1(\Sigma)$.
The isotopy class of a loop $s$ will be denoted by $[s]$. If $b$ is
a boundary component of the surface $\Sigma$, we usually use $b$ to
denote $[b]$.
Given two isotopy classes
$\alpha$ and $\beta$ in $S(\Sigma)$, let $I(\alpha, \beta)$ be their 
\it geometric intersection number \rm which is min$\{ |a \cap b| | 
a \in \alpha, b \in \beta \}$. 
If $f$ is a function defined on $S(\Sigma)$, we define $f(a) = f([a])$.
In particular, the intersection number $I([a], [b])$ is also denoted
by $I(a, [b]) = I([a], b) = I(a,b)$.
We use $\alpha \perp \beta$ to denote  two 
elements $\alpha, \beta \in S(\Sigma)$ so that 
$I(\alpha , \beta) = 1$. And we use $\alpha \perp_0 \beta$
to denote two  elements $\alpha $ and $\beta$  so that
$I(\alpha, \beta) =2$  and 
their algebraic intersection number is zero. Two elements
$\alpha, \beta$ are called disjoint, denoted by $\alpha
\cap \beta = \emptyset$, if $I(\alpha, \beta)  =0$ and $\alpha \neq \beta$.
If $I(\alpha, \beta) \neq 0$, we
say that $\alpha$  \it intersects \rm  $\beta$.
Two isotopic curves
$a, b$ will be denoted by $a \cong b$.

Let $\hat \bold  Q = \bold Q \cup \{\infty\} = \bold Q P^1$.
Two rational numbers $p/q, p'/q'$ satisfying $pq' - p'q = \pm 1$
will be denoted by $p/q \perp p'/q'$. The relation $(\hat \bold Q, \perp )$
is the so called \it modular relation. \rm 
It is well known from elementary number theory that
one may identify $\hat \bold Q$ with the set of cusps in fig. 3.1
so that two cusps are joint by an edge if and only if the
corresponding rational  numbers $r, r'$ satisfy $r \perp r'$. We say
three elements $(\alpha, \beta, \gamma)$ in $\hat \bold Q $  form a
\it triangle \rm if $\alpha \perp \beta \perp \gamma \perp \alpha$
and four distinct elements $(\alpha, \beta, \gamma; \gamma')$ form
a \it quadrilateral \rm if both $(\alpha, \beta, \gamma)$ and
$(\alpha, \beta, \gamma')$ are triangles (see fig. 3.1).

\midspace{0.1cm}
\centerline{\epsfbox{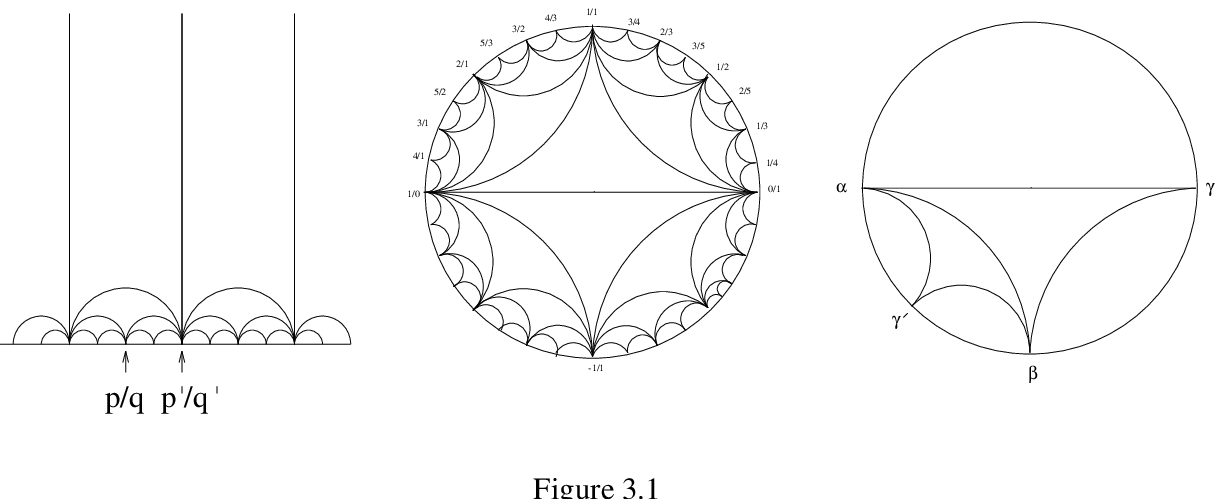}}

We shall always fix an orientation on $\hat \bold Q$ 
so that the triangle
$(0,1, \infty)$ is positively oriented (i.e., the right-hand
orientation in fig. 3.1). A triangle is positively oriented if
it determines the fixed orientation.
The group of orientation preserving bijection of $(\hat \bold Q,
\perp)$ is $PSL(2, \bold Z)$ where the action of the matrices is
given by the fractional linear transformations.

The  importance of $(\hat \bold Q, SL(2, \bold Z))$ in surface
theory was predicted by Grothendieck in [Gr] (page 11, second paragraph).

3.2.
Suppose $\Sigma$ is a  level-1 surface $\Sigma_{1,1}$ or $\Sigma_{0,4}$.
Then there exists a bijection (a slope map) 
$\pi : S'(\Sigma) \to \hat \bold Q = \bold Q \cup \{ \infty 
\}$
so that $\pi(\alpha) = p/q$ and $\pi(\beta) = p'/q'$ satisfy
$pq' - p'q = \pm 1$ if and only if $I(\alpha, \beta)$ = 1 for $\Sigma$ = $\Sigma_{1,1}$ and $I(\alpha, \beta)$ =2
for $\Sigma$ = $\Sigma_{0,4}$. 
This important fact was established by M. Dehn [De] by using Dehn's coding
$ \left(\matrix p \\ q \endmatrix \right)$ 
of classes in $S'(\Sigma)$. 

Here is one way to construct a slope map $\pi : S'(\Sigma) \to \hat \bold  Q$.
It is well known that for the torus
$\Sigma_{1,0}$, $S(\Sigma_{1,0})$ can be naturally identitied with the
set of primitive elements in the first homology group $H_1(\Sigma_{1,1},
\bold Z)$ modulo $\pm 1$. Thus, by fixing a basis for the first homology
group, one constructs a slope map $\pi: S(\Sigma_{1,0}) \to  \hat \bold Q$.
For the 1-holed torus $\Sigma_{1,1}$, let $i$ be an inclusion map
from $\Sigma_{1,1} $ to $ \Sigma_{1,0}$. Then the induced map $i_*$ from
$S'(\Sigma_{1, 1})$ to $S(\Sigma_{1, 0})$ is a bijection preserving
the relation $\perp$. Thus a slope map for $S'(\Sigma_{1,1})$ is
the composition $\pi \circ i_*$. For the 4-holed sphere $\Sigma_{0,4}$,
there exists a natural bijection $P: S'(\Sigma_{0,4}) \to S'(\Sigma_{1,1})$
so that $P(\alpha) \perp P(\beta)$ if and only if $\alpha \perp_0 \beta$.
The bijection $P$ is constructed as follows. Let $T: \Sigma_{1,1}
\to \Sigma_{1,1}$ be a hyperelliptic involution. It is well known that
$T(s) \cong  s$ for any simple loop $s$.
Let
$\Sigma_{0,4}$ be the quotient space $\Sigma_{1,1}/ x \sim T(x)$ 
with a regular
neighborhood of the branch points removed. Then for each $[a]
\in S'(\Sigma_{0,4})$ the inverse image of $a$ in $\Sigma_{1,1}$
consists of two disjoint simple loops $b$ and $T(b)$. Define $P([a]) = [b]$.
Thus a slope map for $\Sigma_{0,4}$ is $\pi \circ i_* \circ P$.

Just like in the modular configuration, for a level-1 surface $\Sigma$, 
we  can talk about triangles and quadrilaterals in $S'(\Sigma)$. 
Furthermore,   when the surface $\Sigma$ is oriented, by making all maps
$\pi$, $i$ and $P$ orientation preserving,  we can talk
about oriented triangles in $S'(\Sigma)$.

\midspace{0.1cm}
\centerline{\epsfbox{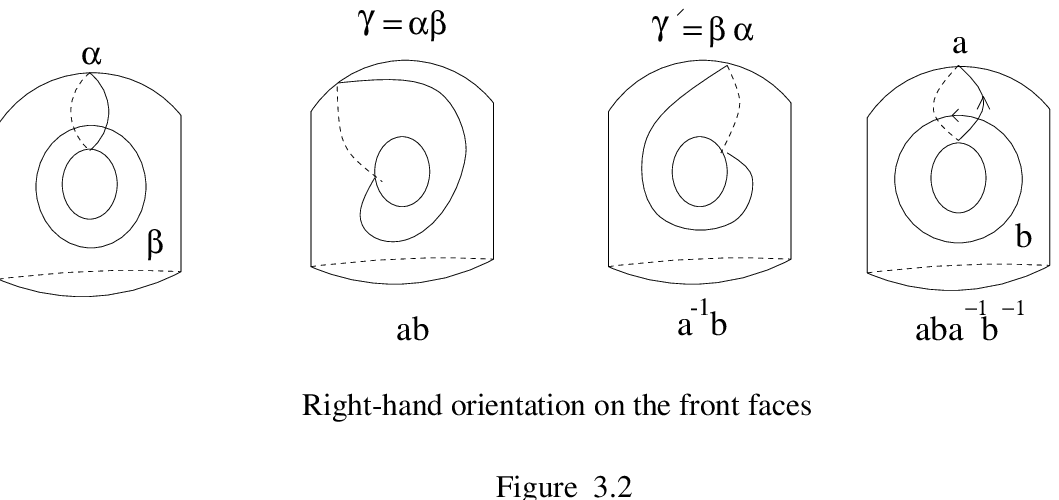}}

The relationship between the fundamental group and the modular
structure  on  $S(\Sigma)$
can be described as follows.
For the 1-holed torus $\Sigma_{1,1}$, if  ($\alpha$, $\beta$, 
$\gamma$; $\gamma')$ is a quadrilateral in $S'(\Sigma_{1,1})$,  then
we can choose generators $a, b$ in the fundamental group $\pi_1(\Sigma_{1,1})$
so that $\alpha$, $\beta$, $\gamma$, and $\gamma'$ 
are represented by $a, b, ab$ and $a^{-1}b$ respectively (see fig. 3.2).
By proposition 2.2, 
the values of an $SL(2,K)$ character
 on $S(\Sigma_{1,1})$ is
determined by its restriction on a triangle.

For the 4-holed sphere $\Sigma_{0,4}$, 
if ($\alpha$, $\beta$, $\gamma$; $\gamma' )$
is a  quadrilateral in $S'(\Sigma_{0,4})$,  then
we can choose three generators $x_1, x_2, x_3$ in the fundamental group
$\pi_1(\Sigma_{1,1})$ so that (1) the four boundary components of the
surface are homotopic to $x_1, x_2, x_3$ and $x_1 x_2 x_3$ and (2) the
classes $\alpha$, $\beta$, $\gamma$ and $\gamma'$ are represented by 
$x_1 x_2$, $x_2 x_3$, $x_1 x_3$ and $x_1 x_2 x_3 x_2^{-1}$ (see fig. 3.3(a)). 
By proposition 2.2 for 3-generator groups, it follows that
an $SL(2,K)$ character defined on $S(\Sigma_{0,4})$ is
determined by its restriction to a triangle and the four boundary
components.

Thus triangles, quadrilaterals and boundary components
in  $S(\Sigma_{0,4})$ and $S(\Sigma_{1,1})$ are exactly 
the elements appeared in lemma 2.2.  
One advantage of using the modular configuration is the
symmetry in the modular configuration. 
For instance, each triangle in the modular configuration $S(\Sigma_{0,4})$
 is invariant under all permutations of the four boundary components
(see fig. 3.3(b)).
As a consequence, there exists a  24-fold
symmetry in the equation (c) in lemma 2.2.

\midspace{0.1cm}
\centerline{\epsfbox{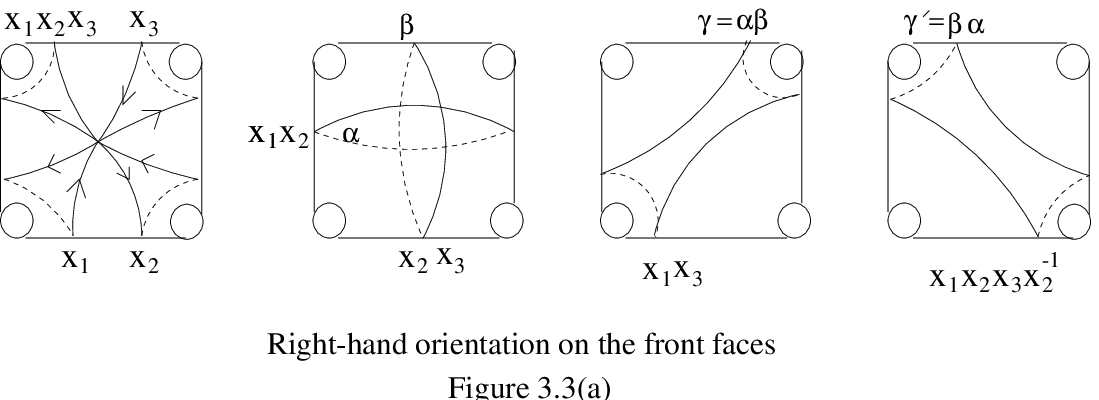}}

\centerline{\epsfbox{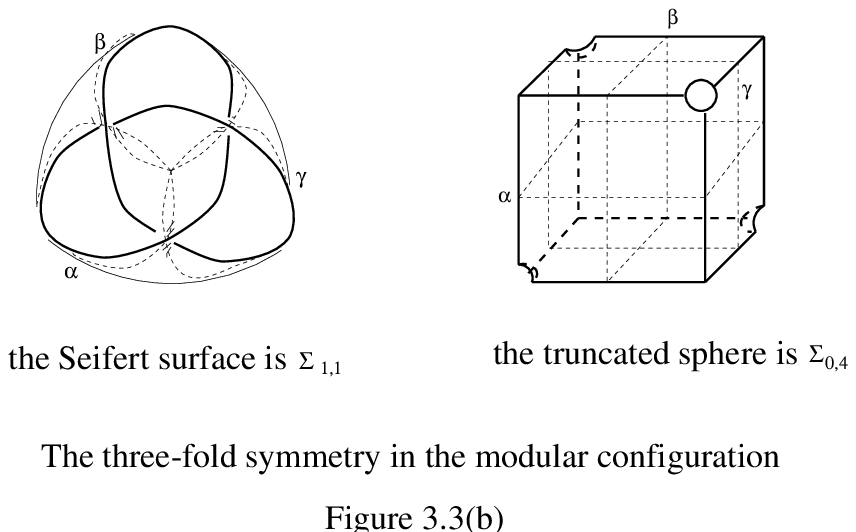}}

In the following,  we shall give a necessary and sufficient condition for a $K$-valued function 
defined on $S(\Sigma)$ to be a character by translating information on the fundamental group 
$\pi_1(\Sigma)$ to $S(\Sigma)$. These results are  certainly well known (see for instance [Go], [Ho
], [Ma], [GM] and others). The only novelty is that it is formulated in terms of the modular configuration.

3.3. For the level-0 surface  $\Sigma_{0,3}$, its fundamental
group is the free group on 2-generator $<x, y>$  where
$x, y$ and $xy$ represent the three boundary components
$b_1, b_2$ and $b_3$. Furthermore, we have $S(\Sigma_{0,3})
=\{ b_1, b_2, b_3\}$.
By proposition 2.2, lemmas 2.3 and 2.5, one obtains the following
result (see [Go]).


{\bf Proposition }.   \it Suppose $\partial \Sigma_{0,3}
= b_1 \cup b_2 \cup b_3$.   Any function $f: S(\Sigma_{0,3}) \to K$ is 
an $SL(2,K)$ character.  Furthermore,
the character is reducible if and only if 
$\sum_{i=1}^3 f^2(b_i) - f(b_1) f(b_2) f(b_3) = 4$.
\rm

By a simple calculation and  by lemma 2.6, one obtains the following.

{\bf Corollary.} \it Under the same assumption as above,

(a) if $f: S(\Sigma_{0,3}) \to K$ satisfies $f^2(b_1) =  4$,
then $f$ is reducible if and only if $f(b_3)  =  f(b_1) f(b_2)/2$
when the characteristic of $K$ is not 2 and $f(b_3) = f(b_2)$ when
the characteristic of $K$ is 2. In particular,  if the
characteristic of $K$ is not 2, and $f^2(b_1) = f^2(b_2) = 4$,
then $f$ is reducible if and only if $f(b_1) f(b_2) f(b_3) = 8$.

(b) If $f: S(\Sigma_{0,3}) \to K$ is a reducible representation
so that $f^2(b_1) \neq 4$, then there exist exactly two $SL(2,K)$
conjugacy classes of $SL(2,K)$ representations whose characters
are $f$. \rm

3.4. For the 1-holed torus, we have,

{\bf Proposition.}  \it Let $b =\partial \Sigma_{1,1}$. A  function 
$f: S(\Sigma_{1,1}) \to K$  is an $SL(2, K)$ character
if and only if the following hold.

$$ \sum_{i=1}^3 f^2(\alpha_i)  - \prod_{i=1}^3 f(\alpha_i) - f(b) = 2 \quad
\text{and} \tag a$$
$$ f(\alpha_3) + f(\alpha_3') =  f(\alpha_1) f(\alpha_2) \tag b$$
where $(\alpha_1, \alpha_2, \alpha_3)$ and $(\alpha_1, \alpha_2, \alpha_3')$
are distinct triangles in $ S'(\Sigma_{1,1})$.

The character $f$ is reducible if and  only if 
$ \sum_{i=1}^3 f^2(\alpha_i) - \prod_{i=1}^3 f(\alpha_i) - 4 = 0$. 
\rm

\it Proof. \rm
The necessity of the conditions 
follow from the trace identities in lemma 2.2 and the choice of
generators for the fundamental group in \S3.2 (see fig. 3.2). 
Due to the modular
relation, if $f_1$ and $f_2$ are two  functions satisfying equation (b) in the
proposition so that they coincide when restricted to
a  triangle,  then $f_1 = f_2$.
Thus the sufficiency of the condition  follows from lemma 2.3.
$\square$

{\bf Corollary}. \it
Suppose $f$ is an irreducible $SL(2,K)$ character defined on the set
$S(\Sigma_{1,1})$. Then either there exists a 3-holed sphere $\Sigma'$
in $\Sigma_{1,1}$ so that the restriction $f|_{ S(\Sigma')}$ is
irreducible or  the characteristic of  $K$ is not 2 and
$f(\partial \Sigma_{1,1}) = -2$ and $f(\alpha) = 0
$ for all $\alpha \in S'(\Sigma_{1,1})$. Furthermore,
if the characteristic of $K$ is not 2, there exists an irreducible
representation $\rho$ of $\pi_1(\Sigma_{1,1})$ so that $tr\rho(\alpha) = 0$
and $tr \rho (\partial \Sigma_{1,1}) = -2$ for all $\alpha \in S'(\Sigma_{1,1})$.
 
\rm

\it Proof. \rm  Let the boundary of $\Sigma_{1,1}$ be $b$.
By the irreducible assumption,
$f(b) \neq 2$. Suppose otherwise that the restriction
of $f$ to each 3-holed sphere is reducible. Since each essential 3-holed sphere
is bounded by $b$ and two copies of  $\alpha \in S'(\Sigma_{1,1})$, by 
lemma 2.5,
we have $ 2f^2(\alpha) + f^2(b) - f^2(\alpha) f(b) = 4$ for all $\alpha \in
 S'(\Sigma_{1,1})$. Since $f(b) \neq 2$, this shows that $f^2(\alpha)  = f(b) + 2$  for
all $\alpha \in  S'(\Sigma_{1,1})$. Now take three elements
$\alpha_i \in  S(\Sigma_{1,1})$ forming a triangle in the modular
configuration. By the above proposition and $f^2(\alpha_i) 
= f(b) + 2$, we obtain
$f^2(\alpha_i) ( 3- f(\alpha_i))  = f^2(\alpha_i)$. Thus either (1) $f(\alpha_i) = 0$ and
$f(b) = -2$ or (2) $f(\alpha_i) = 2$ and $f(b) = 2$. But the case (2)
and the case (1) when the characteristic of $K$ is 2 are
excluded by the irreducibility assumption.  

If the characteristic of $K$ is not 2, then one constructs an 
$SL(2,K)$ representation of $\pi_1(\Sigma_{1,1})$  satisfying
the above condition as follows. First we note that
the unit quarternion group 
 of eight elements $\{\pm 1, \pm i,$$ \pm j,$$ \pm k\}$ is a subgroup of
$SL(2,K)$ where $i = \left(\matrix 0& -1\\1 & 0 \endmatrix \right)$
and $j = \left(\matrix \sqrt{-1} & 0 \\ 0 & - \sqrt{-1} \endmatrix \right)
$. Now the character of any representation of $\pi_1 (\Sigma_{1,1})$
onto the unit quarternion group 
takes value zero on $ S'(\Sigma_{1,1})$ and $-2$ on the
boundary component.
Any two such representations are $SL(2,K)$ conjugate.
$\square$

3.5. For the 4-holed torus, we have,

{\bf Proposition.} \it Let $\partial \Sigma_{0,4} = \cup_{i=1}^4 b_i$. A
function $f: S(\Sigma_{0,4}) \to K$ is an $SL(2, \bold K)$ 
character if and only if for each  triangle $(\alpha_1,
\alpha_2, \alpha_3)$ in $S'(\Sigma_{0,4})$ the following hold.

$$\sum_{i=1}^3 f^2(\alpha_i) +
 \prod_{i=1}^3 f(\alpha_i) + \sum_{r=1}^4 f^2(b_r) + \prod_{r=1}^4 f(b_r) -
\sum_{(i,r,s) \in P}f(\alpha_i)f(b_r)f(b_s) -4 = 0 \tag a $$
$\text{where $P$ =\{ $(i,r,s)$ $|(\alpha_i, b_r, b_s)$ bounds
a $\Sigma_{0,3}$\} } $ and, 
$$ f(\alpha_3) + f(\alpha_3') = -f(\alpha_1) f(\alpha_2)
 + f(b_i)f(b_j)  + f(b_k)f(b_l) \tag b $$
where $(3,i,j)$ and $(3,k,l)$ are in $P$.

A character is irreducible if and only if there is a 3-holed sphere
so that the restriction of $f$ to the 3-holed sphere is irreducible.
\rm

\it Proof. \rm The necessity of the conditions 
follow from the trace identities in lemma 2.2 and the choice of
generators for the fundamental group in \S3.2 (see fig. 3.3). 

To show the sufficiency of the conditions, we first
note that  by the modular relation and the iteration equation (b),
each $f$ is determined by its restriction to
a 7-element set $\{\alpha_i, b_r| i=1,2,3; r=1,2,3,4\}$
where $\alpha_i$'s form a triangle. We choose a set of
generators $x_i$ for the fundamental group $\pi_1(\Sigma_{0,4})$
as in \S3.2 so that $\alpha_k$ is represented by $x_i x_j$,
$i \neq j \neq k \neq i$.  Let $t_i = f(b_i)$ 
and $t_{ij} = f(\alpha_k)$ where
$i=1,2,3$, $(i,j) = (1,2), (2,3), (3,1)$ and $k \neq i \neq j \neq k$.
By lemma 2.3, we find $A_i \in SL(2,K)$ so that $tr(A_i) = t_i$ and
$tr(A_{ij}) = t_{ij}$. By equation (a), $f(b_4)$ is a root of
the quadratic equation $x^2 -Px + Q = 0$ where $P$ and $Q$ are the
same as in lemma 2.2. By lemma 2.2 (c), we may assume after change
$(A_1, A_2, A_3)$ to $(A_1^{-1}, A_2^{-1}, A_3^{-1})$ if necessary
that $f(b_4)$ = $tr(A_1A_2A_3)$. Define a representation of
$\pi_1(\Sigma_{0,4})$ by sending the generator $x_i$ to $A_i$.
Then the character of this representation
and the function $f$ take the same values on the seven specific
elements. Thus they are the same.

The irreducibility condition follows from corollary 2.5(a).
$\square$.

{\bf Corollary.}  \it 
(a)
Suppose $f$ is an irreducible character defined on
$S(\Sigma_{0,4})$.
Then for any element $\alpha \in S'(\Sigma_{0,4})$,
there exists an element $\beta$ with
$I(\alpha, \beta) =2$ and a 3-holed sphere $\Sigma'$ bounded by
$\beta$ and two components of $\partial \Sigma$ so that the
restriction $f|_{ S(\Sigma')}$ is irreducible.

(b) 
Let $b$ be  a boundary component of $\partial \Sigma_{0,4}$ and
$\rho$ be an irreducible representation of $\pi_1(\Sigma_{0,4})$.
If
$\rho(b) \neq \pm id$, in particular if $tr \rho(b) \neq \pm 2$,
then there exists a level-0 subsurface $\Sigma'$ having $b$ as
a boundary component so that the restriction of $\rho$ to
$\Sigma'$ is irreducible. 

(c)  If an $SL(2,K)$ representation of $\pi_1(\Sigma_{0,n})$
is reducible on each level-0 subsurface group, then the
representation reducible.

(d) If $f$ is a reducible character on $S(\Sigma_{0,4})$, then
for any quadrilateral $(\alpha, \beta, \gamma; \gamma')$ in
$S'(\Sigma_{0,4})$, $f(\gamma) = f(\gamma')$.
\rm

\it Proof. \rm  
To prove (a),  suppose otherwise that $f$ is reducible on each level-0
subsurface $\Sigma'$ in $\Sigma_{0,4}$ bounded by a simple loop intersecting 
$\alpha$ at two points. Then take a quadrilateral with vertices
$(\alpha, \beta, \alpha \beta; \beta \alpha)$
Choose three generators $\{x_1, x_2, x_3\}$ for $\pi_1(\Sigma_{0,4})$ 
as in fig. 3.3
so that $\beta$, $\alpha \beta$, $\beta \alpha$
are represented by $x_2x_3$, $x_1 x_3$, and
$x_1 x_2 x_3 x_2^{-1}$.
Let the $SL(2,K)$ representation corresponding to the character $f$
send $x_i$ to the matrix $A_i$ and let $A_4 = A_1 A_2 A_3$. 
Since the representation is reducible over six 3-holed spheres
bounded by $\beta, \alpha \beta$ and $\beta \alpha$, we
obtain the following reducible equations by lemma 2.5. Namely,
$\Delta (A_1, A_3) = \Delta (A_2, A_3) = \Delta (A_2, A_1A_3) =
 \Delta (A_3, A_1A_2A_3A_2^{-1}) = 0$.

It turns out
that these four  equations imply that the group $<A_1, A_2, A_3>$
is reducible. Indeed, if $A_3 = \pm id$, 
then $\Delta (A_2, A_1A_3) = \Delta (A_2, A_1 )
=0$ shows that the group $<A_1, A_2, A_3> =<A_1, A_2, \pm id>$ is
reducible by lemma 2.5. Suppose $A_3$ has exactly one eigenspace. Then
$\Delta (A_1, A_3)$ $=$ $\Delta (A_2, A_3)$ $=0$ 
imply that the eigenspace is fixed
by  both $A_1$ and $A_2$. Thus the group is again reducible. Finally,
suppose $A_3$ has exactly two distinct eigenspaces $L_1$ and $L_2$ so
that none of $L_i$ is fixed by both $A_1$ and $A_2$.
Since $\Delta(A_3, A_i) = 0$ for $i=1,2$, we may assume that
$A_i(L_i) = L_i$ for $i=1,2$. By $\Delta (A_3, A_1A_2A_3A_2^{-1}) =0$,
one of the eigenspace $L_i$ is fixed by $ A_1A_2A_3A_2^{-1}$. If
$ A_1A_2A_3A_2^{-1}(L_2) = L_2$, then $A_1(L_2) =L_2$. Thus the group
has a common eigenspace $L_2$ and is reducible. If
 $A_1A_2A_3A_2^{-1}(L_1) = L_1$, then  $A_3(A_2^{-1} (L_1)) = A_2^{-1}(L_1)$.
Thus either $A_2^{-1}(L_1) = L_1$ or $A_2^{-1}(L_1) = L_2$. In the
first case, $L_1$ is a common eigenspace for the group. The second case
implies $L_1 = L_2$ which is absurd. In summary, we have shown that
the group is reducible which contradicts the assumption.

To prove (b), suppose otherwise that the restrictions of $\rho$ to
all level-0 subsurfaces having $b$ as a boundary component are reducible.
Let us choose a set of generators $x_1, x_2, x_3$ for $\pi_1(\Sigma_{0,4})$
as in \S3.2 so that $x_1$ corresponds to $b$. Let $\rho(x_i)$
be the matrix $A_i$ in $SL(2, K)$. Then
by the assumption $\Delta (A_1, A_2) =\Delta (A_1, A_3) = \Delta (A_1, A_1A_2A_3) = 0$.
Since $\rho(b) \neq \pm id$, the matrix $A_1$ has at most two eigenspaces.
If it has exactly one eigenspace, then $\Delta (A_1, A_2) = \Delta (A_1, A_3) =0$
implies that the  eigenspace is invariant under both $A_2$ and $A_3$.
This contradicts the irreducible  assumption. If $A_1$ has two
distinct eigenspaces $L_2$ and $L_3$ so that $A_i(L_i) = L_i$ for
$i=2,3$, then due to $\Delta (A_1, A_1 A_2 A_3) =0$,
 one of $L_i$ is invariant under $A_1A_2A_3$.
But this again implies that $L_i$ is
a common eigenspace of $A_1$, $A_2$ and $A_3$.

To prove (c), we first note that the result holds for $n=3$
by part (a).  For $n>3$,
take generators $x_1, ..., x_{n-1}$ for the
free group $\pi_1(\Sigma_{0,n})$ so that the boundary components
are freely homotopic to $x_i$ or $x_1 ... x_{n-1}$ as in figure 2.1. 
 Now each
3-generator subgroup $<x_i, x_j, x_k>$ lies in a level-1 subsurface
subgroup. Thus the restriction of the representation to the
3-generator subgroup is reducible. 
By corollary 2.5(b), this shows that the representation is reducible.

To see part (d), we may assume that $f$ is the character of
a diagonalizable representation $\rho$. Choose a set of generators
$\{x_1, x_2, x_3\}$ for the fundamental group so that $\gamma$ and
$\gamma'$ are represented by $x_1x_3$ and $x_1x_2x_3x_2^{-1}$
as in fig. 3.3. Then $f(\gamma) = tr \rho(x_1x_3) = tr \rho(x_1x_2x_3x_2^{-1})
= f(\gamma')$.
$\square$

3.6. In this section, we give a different
interpretation of the triangles and quadrilaterals in the set 
$S(\Sigma)$. This new interpretation is the basis for us 
in dealing with simple loops in the rest of the paper.

We begin by introducing  some notations.
Recall that surfaces are oriented.
If $a$ and $b$ are two arcs intersecting  transversely at a point $p$,
then the \it resolution of \rm  $a \cup b$  \it at $p$ from $a$ to $b$ \rm
is defined as  follows. Fix any orientation  on $a$ and use the orientation
on the surface to determine an orientation on $b$. Then resolve
the intersection according to the orientations (see fig. 3.4).
The resolution is
evidently independent of the choice of the orientations on $a$. If
$\alpha \perp \beta$ or $\alpha \perp_0 \beta$, take $a \in \alpha$,
$b \in \beta$ so that $|a \cap b| = I(\alpha , \beta)$. Then the
curve obtained by resolving all intersection points in
$a \cap b$ from $a$ to $b$ is again a simple loop denoted by $ab$.
We define $\alpha \beta$ to be the isotopy class of $ab$. It follows
from the definition that when $\alpha \perp \beta$ then $\alpha \beta
\perp \alpha, \beta$ and when $\alpha \perp_0 \beta$ then $\alpha
\beta \perp_0 \alpha, \beta$. In particular, if $\Sigma$ has level-1,
then the positively
oriented triangles in $S'(\Sigma)$ are $(\alpha, \alpha \beta, \beta)$
where $\alpha \perp \beta$ or $\alpha_0 \perp \beta$. Also the quadrilaterals
are $(\alpha, \beta, \alpha \beta; \beta \alpha)$.
Let $N(a)$ and $N(b)$ be two small regular neighborhoods of $a$ and $b$. Then
$N(a \cup b) = N(a) \cup N(b)$ is homeomorphic to $\Sigma_{1,1}$ when
$\alpha \perp \beta$ and to $\Sigma_{0,4}$ when $\alpha \perp_0 \beta$.
We use $\partial(\alpha, \beta)$ to  denote the set  of isotopy classes
of the curves in $\partial N(a \cup b)$.

\midspace{0.1cm}
\centerline{\epsfbox{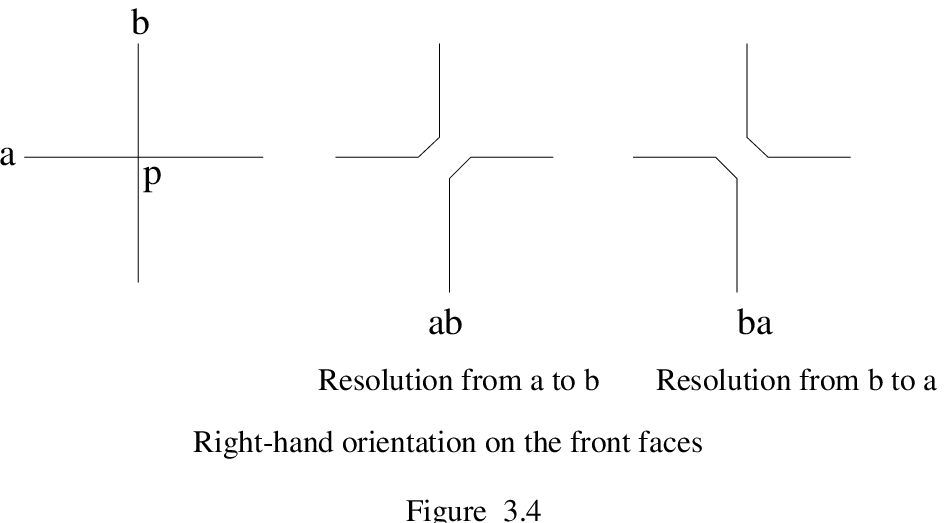}}

In terms of these notations, the equations (b) in
propositions 3.4 and 3.5 say that if $\alpha \perp \beta$, then
$f(\alpha \beta)$ is determined by $f(\alpha), f(\beta)$ and
$f(\beta \alpha)$, and if $\alpha \perp_0 \beta$, then
$f(\alpha \beta)$ is determined by the values of $f$ on
$\{\alpha, \beta, \beta \alpha\}$ and $\partial(\alpha, \beta)$.
More precisely, a function $f: S(\Sigma) \to K$ is a
\it trace function \rm if
and only if  it satisfies:

(1) $f(\alpha \beta) + f(\beta \alpha) = f(\alpha) f(\beta)$,
for
$\alpha \perp \beta$,

(2) $f^2(\alpha) + f^2(\beta) + f^2(\alpha \beta) - f(\alpha) f(\beta)
f(\alpha \beta) = f(\partial (\alpha, \beta)) + 2$ for $\alpha \perp \beta$,

(3) $f^2(\alpha) + f^2(\beta) + f^2(\alpha \beta) + f(\alpha)f(\beta) f(\alpha
\beta) + \Pi_{i=1}^4 f(\gamma_i) + \Sigma_{i=1}^4f^2(\gamma_i) -$
$f(\alpha)(f(\gamma_1)f(\gamma_2) + f(\gamma_3)f(\gamma_4))$
$-f(\beta)(f(\gamma_1)f(\gamma_3) + f(\gamma_2) f(\gamma_4))$
$- f(\alpha \beta)(f(\gamma_1)  f(\gamma_4) + f(\gamma_2)f(\gamma_3)) - 4 = 0$,
for $\alpha \perp_0 \beta$, $\partial(\alpha, \beta) = \{ \gamma_1, ...,
\gamma_4\}$ so that $(\alpha, \gamma_1, \gamma_2)$ and $(\beta, \gamma_1,
\gamma_3)$ bound level-0 subsurfaces,

(4) $f(\alpha \beta) + f(\beta \alpha) = -f(\alpha) f(\beta)
+ f(\gamma_1) f(\gamma_4) + f(\gamma_2) f(\gamma_3)$ under the
same assumption as in (3).

3.7.
One of the main reduction lemma for simple loops is 
the following
which generalizes Lickorish's lemma 2 in [Li1]. See [Lu2] lemma 7
for a proof.

{\bf Lemma.} \it Suppose $\gamma_1, ..., \gamma_m$ are pairwise
disjoint classes in $S(\Sigma)$. If $\alpha \in S(\Sigma)$
intersects $\gamma_1$ and is
not $\perp$ or $\perp_0$ related to $\gamma_1$, 
then 
$\alpha = \beta_1 \beta_2$ with $\beta_1 \perp \beta_2$ or
$\beta_1 \perp_0 \beta_2$
so that 

(1) $I(\beta_i, \gamma_1) < I(\alpha, \gamma_1)$, 
$I(\beta_i, \gamma_j) \leq I(\alpha, \gamma_j)$, 
$I(\beta_2 \beta_1, \gamma_1) < I(\alpha, \gamma_1)$ and
$I(\beta_2 \beta_1, \gamma_j) \leq I(\alpha, \gamma_j)$ for
all $i=1,2$ and $j =2, ..., m$, and,

(2) if $\beta_1 \perp_0 \beta_2$, then for each element $\delta
\in \partial (\beta_1, \beta_2)$, $I(\delta, \gamma_1) < I(\alpha,
\gamma_1)$ and $I(\delta, \gamma_j) \leq I(\alpha, \gamma_j)$
for $j = 2,..., m$.

\rm 

The lemma says that one can ``simplify" $\alpha$ unless
$\alpha \cap \gamma_1 = \emptyset$, or $\alpha \perp \gamma_1$,
or $\alpha \perp_0 \gamma_1$. In particular, if we set
$\Cal G_0 =\{ \alpha \in S(\Sigma) |$ for each $i$,  either
$\alpha \cap \gamma_i = \emptyset$, or $\alpha \perp \gamma_i$,
or $\alpha \perp_0 \gamma_i$\}, then by induction on
$(I(\alpha, \gamma_1), ..., I(\alpha, \gamma_m))$, we have
 $S(\Sigma) = \cup _{n=0}^{\infty} \Cal G_n$
where $\Cal G_{n+1} = \Cal G_n \cup \{ \alpha |$ 
$\alpha = \beta_1 \beta_2$ either (1) $\beta_1 \perp \beta_2$ and
$\{ \beta_1, \beta_2, \beta_2 \beta_1\} \subset \Cal G_n$ or (2)
$\beta_1 \perp_0 \beta_2$ and  \{$\beta_1, \beta_2, \beta_2 \beta_1
\} \cup \partial(\beta_1, \beta_2) \subset \Cal G_n$\}.

By the remark in the last two paragraphs in \S3.6, we obtain,

{\bf Corollary.} \it (a) Let $\gamma_1, ..., \gamma_m$ be  pairwise
disjoint classes in $S(\Sigma)$ and 
$f$ and $g$  be two trace functions
on $S(\Sigma)$. If
$f(\alpha) = g(\alpha)$ for each class $\alpha \in S(\Sigma)$
so that for each $i$ either $\alpha \cap \gamma_i = \emptyset$,
or $\alpha \perp \gamma_i$ or $\alpha \perp_0 \gamma_i$, 
then $f = g$. 

(b) Let $\gamma_1$ and $\gamma_2$ are two disjoint classes in
$S'(\Sigma)$ so that $\gamma_1$ bounds a $\Sigma_{1,1}$ 
and $\gamma_2$ lies in $\Sigma_{1,1}$. If $f$ and $g$ are
two trace functions on $S(\Sigma)$ so that $f(\alpha) = g(\alpha)$
for all $\alpha \perp \gamma_2$ and $\alpha \perp_0 \gamma_1$,
then $f = g$.

\rm

Indeed, in part (b), if $\alpha \perp_0 \gamma_1$, then $\alpha \perp_0
\gamma_2$ cannot occur. Thus part (b) follows from part (a).

3.8. As a second consequence of the above lemma, we obtain the
following result which will be used in \S6. Part (a) of the
corollary was known to many people [Li2].

{\bf Corollary.} \it 
(a) Given two non-separating classes
$\alpha$ and $\alpha'$, there exists a sequence of non-separating
classes $\{ \alpha_i | i=1,..., m\}$ starting from $\alpha$ and
ending at $\alpha'$ so that 
$\alpha_i \perp \alpha_{i+1}$ for all $i$.

(b) Given two essential 1-holed tori (resp. level-1 surfaces)
 $T$ and $T'$ in $\Sigma$,
there exists a sequence of essential 1-holed tori (resp. level-1
surfaces) \{$T_i$ \} starting
from $T$ and ending at $T'$ so that  $S'(T_{i}) \cap S'(T_{i+1}) \neq \emptyset$
for all $i$.
\rm

\it Proof. \rm  Let us denote two classes $\alpha$ and $\alpha'$
satisfying the conclusion of (a) by $\alpha \sim \alpha'$.  We use
the induction on $I(\alpha, \alpha')$ to prove part (a). 
Clearly if $\alpha \cap \alpha' = \emptyset$ or $\alpha \perp \alpha'$,
then $\alpha \sim \alpha'$. If $\alpha \perp_0 \alpha'$, since
both $\alpha$ and $\alpha'$ are non-separating, one of the
element $\beta \in \partial (\alpha, \alpha')$ is non-separating.
Thus $\alpha \sim \beta \sim \alpha'$. In the remaining cases,
by lemma 3.7, we can write $\alpha = \beta_1 \beta_2$ where either
$\beta_1 \perp \beta_2$ or $\beta_1 \perp_0 \beta_2$ and
$I(\beta_i, \alpha') < I(\beta_i, \alpha)$, $i=1,2$. Since $\alpha$
is non-separating, one of $\beta_1$ or $\beta_2$, say $\beta_1$, is
again non-separating. Thus by the induction hypothesis, $\beta_1 \sim 
\alpha'$. But $\beta_1 \perp \alpha$ or $\beta_1 \perp_0 \alpha$.
Thus $\alpha \sim \alpha'$.

To see part (b), take $\alpha \in S'(T)$ and $\alpha' \in S'(T')$.
Let $T_1 = T$, $T_{m+1} = T'$ and $T_i$ be the 1-holed torus containing
both $\alpha_i$ and $\alpha_{i+1}$. Then the result follows. The
result for level-1 surfaces $T$ and $T'$ is simpler. We omit the
proof.

3.9. It is shown in section 3 of [Lu1]  that there exists a finite set
$F_0 \in S(\Sigma)$ so that $S(\Sigma) = \cup_{n=0}^{\infty} F_n$
where $F_{n+1} = F_n \cup \{ \alpha |$
$\alpha = \beta_1 \beta_2$ either (1) $\beta_1 \perp \beta_2$ and
$\{\beta_1, \beta_2, \beta_2 \beta_1\} \subset  F_n$ or (2)
$\beta_1 \perp_0 \beta_2$ and  \{$\beta_1, \beta_2, \beta_2 \beta_1
\} \cup \partial(\beta_1, \beta_2) \subset F_n$\}.
In particular, if $f$ is a trace function defined on $S(\Sigma)$,
then $f$ is  algebraically determined by $f|_F$. This shows 
the following result analogous to proposition 2.2.

{\bf Propositoin.} \it There exists a finite set of isotopy classes of
simple loops in each compact orientable surface  $\Sigma$
so that  $SL(2,K)$
characters  and trace functions 
on $S(\Sigma)$ are algebraically determined by the
restrictions of the characters on the finite set. \rm


\bigskip

{\bf \S4.  $SL(2,K)$ Characters on  the 5-holed Sphere}

4.1.
We shall use the following terminologies. If $f$ is a 
trace function on $S(\Sigma)$ and $\Sigma'$ is a subsurface of  $\Sigma$, then we 
call
$f|_{S(\Sigma')}$ the \it restriction \rm of $f$ to the subsurface
$\Sigma'$. If $f|_{S(\Sigma')}$ is irreducible, we say $\Sigma'$ is an
\it irreducible \rm subsurface with respect to $f$. We say a subsurface
$\Sigma'$ is \it bounded \rm by $\alpha_1, ..., \alpha_k \in$ $S'(\Sigma)$ if 
$\partial \Sigma' = a_1 \cup ...\cup a_k \cup b_1 \cup ...\cup b_m$ 
so that $a_i \in \alpha_i$ and $b_j \subset \partial \Sigma$.

The goal of this section is to prove the  following case of theorem 1.2
for the 5-holed sphere.

{\bf Theorem.} \it Suppose $K$ is a quadratically closed field. 
If $f$ is a $K$-valued trace function on $S(\Sigma_{0,5})$ so that
either $f$ is reducible on all level-0 subsurfaces or $f$
is irreducible on a level-0 subsurface bounded by two disjoint
elements in $S'(\Sigma_{0,5})$, then $f$ is an $SL(2,K)$ character. 
\rm

In \S5.2, we prove that if $f$ is a trace function on $S(\Sigma_{0,5})$
which does not satisfy the conditions in the above theorem, then
$f$ is exceptional.

4.2. \it The Pentagon Relations \rm

Given five pairwise distinct elements $\alpha_1$, ..., $\alpha_5$ in  $S'(\Sigma_{0,5})$
so that
$\alpha_i$ $\cap$ $\alpha_j$ $= \emptyset$  for $i \neq j \pm 1$  mod 5, it is shown
in [Lu3] that $\alpha_i \perp_0 \alpha_{i+1}$ for all indices $i$ mod 5
(see fig. 4.1).
We say \{$\alpha_1$, ..., $\alpha_5$\} forms a \it
 pentagon  \rm in this case. If \{$\alpha_1$,...,
$\alpha_5$\} forms a  pentagon, then the following conditions hold.

(a) $(\alpha_i \alpha_j) \alpha_k = \alpha_i (\alpha_j \alpha_k)$.

(b) $\alpha_i \alpha_{i+1} \alpha_{i+2} = \alpha_{i+3}\alpha_{i+4}$.

(c) $(\alpha_i \alpha_j) \cap (\alpha_i \alpha_k) = \emptyset$
and $(\alpha_j \alpha_i) \cap (\alpha_k \alpha_i) = \emptyset$,
$i \neq j \neq k \neq i$.

(d) $\alpha_i \alpha_j \alpha_k = \alpha_j \alpha_i \alpha_k$ and
 $\alpha_k \alpha_j \alpha_i = \alpha_k \alpha_i \alpha_j$,
 if $i \neq j \pm 1 $ mod 5.

(e) $\alpha_i \alpha_j \alpha_i = \alpha_j$ if $i = j \pm 1$ mod 5.

These can be verified easily using the definition of resolution or see [Lu2]
or [Lu3] for a proof.
Note that, by (c), if
 $\{\alpha_1, \alpha_2, \alpha_3, \alpha_4, \alpha_5\}$ 
forms a pentagon, then \{$ \alpha_1 \alpha_2,$$
\alpha_2,$$ \alpha_3 \alpha_2,$$ \alpha_4,$$ \alpha_5\}$
is also a pentagon (see fig. 4.1).

4.3.
In this section, we prove that if $f$ is a trace function which is
irreducible on a 3-holed sphere bounded by two elements
$\alpha_2$ and $\alpha_5$ in $S'(\Sigma_{0,5})$, then $f$ is a
character.

Let $X$ and $Y$ be the level-1 subsurfaces bounded by $\alpha_2$ and $\alpha_5$ respectively.
By proposition 3.5, we find two representations $\rho_X$ and $\rho_Y$
of $\pi_1(X)$ and $\pi_1(Y)$ respectively so that $\chi_{\rho_X} = f|_{S(X)}$
and $\chi_{\rho_Y} = f|_{S(Y)}$. The restrictions of $\rho_X$ and
$\rho_Y$ to $\pi_1(X \cap Y)$ have the same character by the construction 
and both are irreducible.
Thus by lemma 2.4,  these two restrictions are conjugate. After conjugate
$\rho_X$, we may assume that $\rho_X |_{\pi_1(X \cap Y)} = \rho_Y |_{\pi_1(
X \cap Y)}$. This defines a representation $\rho: \pi_1(\Sigma_{0,5})
\to SL(2, K)$ so that its restrictions to $\pi_1(X)$ and $\pi_1(Y)$
are $\rho_X$ and $\rho_Y$ respectively. Let $g$ be the character of $\rho$.
Then $f(\alpha) = g(\alpha)$ for all $\alpha \in S(X) \cup S(Y)$. The 
goal is to show that $f=g$ using the irreducibility condition.

By corollary 3.7 applied to $f$ and $g$, it suffices to prove that for
each $\alpha_1$ so that $\alpha_1$$\perp_0 \alpha_2$ and $\alpha_1 \perp_0 \alpha_5$, we have
$f(\alpha_1) = g(\alpha_1)$. 
To this end, we extend $\{\alpha_1, \alpha_2, \alpha_5\}$ to a set
\{$\alpha_1$, ..., $\alpha_5$\} forming a pentagon by setting $\alpha_3 
= \partial(\alpha_1, \alpha_5) \cap S'(\Sigma_{0,5})$ and
$\alpha_4 = \partial(\alpha_1, \alpha_2) \cap S'(\Sigma_{0,5})$.
We shall use proposition 3.5 to derive a system of linear
equations and show that the system has a unique solution and that both
$f(\alpha_1)$ and $g(\alpha_1)$ are solutions.

We begin by introducing some notations. Let $h$ be a trace function
defined on $S(\Sigma_{0,5})$. Given a set of indices $i_1,..., i_k$, $k=1,2,3$,
let $x_{i_1...i_k} = h(\alpha_{i_1}...\alpha_{i_k}) $ when
$\alpha_{i_1} ... \alpha_{i_k}$ is not in $S(X) \cup S(Y)$ and
let $a_{i_1....i_k} = h(\alpha_{i_1} ... \alpha_{i_k})$ when
$\alpha_{i_1} ... \alpha_{i_k} \in S(X) \cup S(Y)$.  Let $\beta_i$
be the component of $\partial \Sigma_{0,5}$ so that $\{\alpha_{i-1},
\beta_{i}, \alpha_{i+1}\}$ bounds a 3-holed sphere (indices mod 5)
and let $b_i = h(\beta_i)$  (see fig. 4.1). Let $\tau$ be the
orientation preserving involution of $\Sigma_{0,5}$ so that $\tau$ sends
$\alpha_{1+i}$ to $\alpha_{1-i}$ and $\beta_{1+i}$ to $\beta_{1-i}$.

Now we derive equations for $x_i$, $x_{ij}$ and $x_{ijk}$ with coefficient
in  $h(\alpha)$'s where $\alpha \in S(X) \cup S(Y)$
 using proposition 3.5.

\midspace{0.1cm}
\centerline{\epsfbox{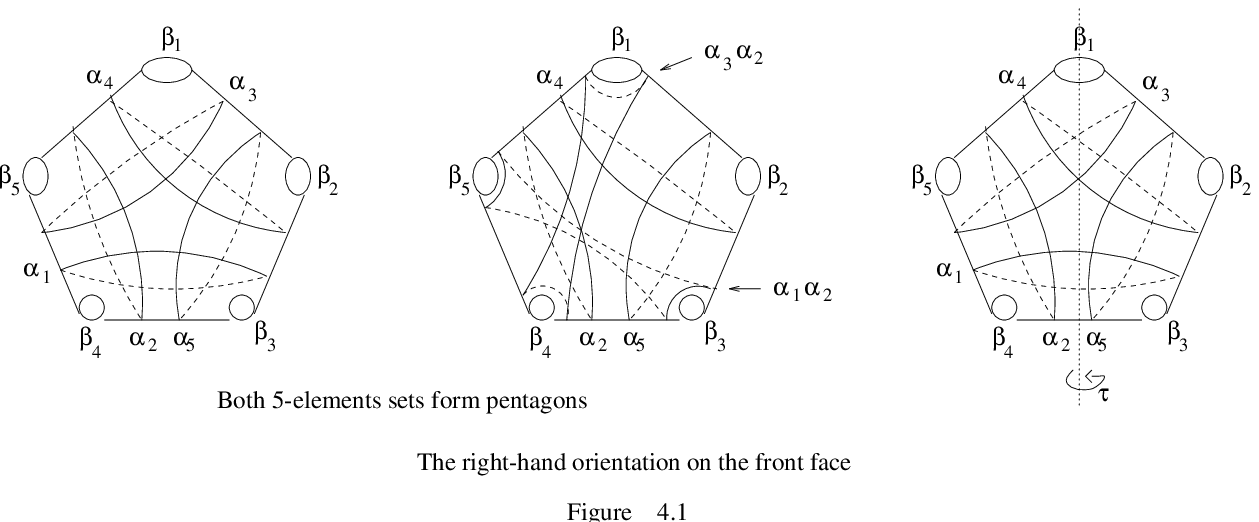}}

Since $\alpha_3 \perp_0 \alpha_4$ and $\partial(\alpha_3, \alpha_4)
\cap S'(\Sigma_{0,5}) = \alpha_1$, by proposition 3.5 (b), we obtain 
$$ h(\alpha_3 \alpha_4) + h(\alpha_4 \alpha_3) = -h(\alpha_3) h(\alpha_4)
+ h(\alpha_1) h(\beta_1) + h(\beta_2)h(\beta_5).$$  This is the same as,
 
$$ x_{34} + x_{43} - b_1 x_1 = p_1. \tag 1$$

Here and below, $p_i$ always denotes some polynomial with integer
coefficient in $h(\alpha)$'s where $\alpha \in S(X) \cup S(Y)$.

Since $\alpha_1 \perp_0 \alpha_2$ with $\partial (\alpha_1, \alpha_2)
\cap S'(\Sigma_{0,5}) = \alpha_4 \in S(X)$, we obtain,

$$ x_{12} + x_{21} + a_2 x_1 = p_2.   \tag 2$$

Apply the involution $\tau$ to the equation (2), we obtain

$$x_{15} + x_{51} + a_5 x_1 = p_3.   \tag 3$$

Since $\alpha_2 \perp_0 \alpha_3 \alpha_4,$ $\alpha_2 ( \alpha_3 \alpha_4)
= \alpha_5 \alpha_1$ and $\partial(\alpha_2, \alpha_3 \alpha_4)
\cap S'(\Sigma_{0,5}) = \alpha_5 \alpha_4 \in S(X)$, we obtain,

$$x_{51} + x_{342} + a_2 x_{34} = p_4.  \tag 4$$

Since $\alpha_4 \perp_0  \alpha_3 \alpha_2,$ $\alpha_4( \alpha_3 \alpha_2)
= \alpha_1 \alpha_5$, $\alpha_3 \alpha_2 \alpha_4 = \alpha_3 \alpha_4 \alpha_2$
and $\partial(\alpha_4, \alpha_3 \alpha_2) \cap S'(\Sigma_{0,5})= \alpha_1 \alpha_2$, we obtain,

$$ x_{15} + x_{342} - b_1x_{12} = p_5. \tag 5$$

Subtracting (4) by (5), we obtain,

$$ x_{51} - x_{15} + a_2 x_{34} + b_1 x_{12} = p_6. \tag 6$$

Apply the involution $\tau$ to equation (6), we obtain,

$$x_{21} - x_{12} + a_5x_{43} + b_1 x_{15} = p_7. \tag 7$$

Since $\alpha_1 \alpha_2 \perp_0 \alpha_5$, $\alpha_5 ( \alpha_1 \alpha_2)
=\alpha_3 \alpha_4$ and $\partial(\alpha_1 \alpha_2, \alpha_5)
\cap S'(\Sigma_{0,5}) = \alpha_3 \alpha_2 \in S(Y)$, we obtain,

$$x_{34} + x_{125} + a_5 x_{12} = p_8. \tag 8$$

Apply the involution $\tau$ to (8)
 and use the fact that $x_{125} = x_{152}$ (due to \S4.1(d)), we
obtain,

$$ x_{43} + x_{125} + a_2 x_{15} = p_9. \tag 9$$

Subtracting (8) by (9) gives,

$$x_{34} - x_{43} - a_2 x_{15} + a_5 x_{12} = p_{10}. \tag 10$$

Now consider the system of linear equations (1), (2), (3), (6), (7), and
(10). By (1), (2) and (3), we obtain $x_{21} = -a_2 x_{1} - x_{12} + p_2$,
$x_{51} = -a_5x_1 - x_{15} + p_3$ and $x_{43} = b_1 x_1 - x_{34} + p_1$. 
Thus, after
substituting these into (6), (7) and (10), we obtain the following
system of linear equations.

$$
\aligned
 b_1 x_{12} - 2 x_{15} + a_2 x_{34} - a_5 x_1 = p_{11}\\
 -2x_{12} + b_1 x_{15} - a_5 x_{34} + (a_5 b_1 -a_2)x_1 = p_{12}\\
 a_5x_{12} - a_2 x_{15} + 2x_{34} - b_1x_1 = p_{13}. 
\endaligned  \tag 11
$$

Let $A$ be the $3 \times 4$ coefficient matrix of the linear system
and $B$ be the $3 \times 3$ submatrix obtained from $A$ by removing
the $4$-th column. Then the determinant of $B$ is $2 \Delta$
where $\Delta = a_2^2 + a_5^2 + b_1^2 - a_2 a_5 b_1 -4$.
Suppose $B^*$ is the adjoint matrix of $B$. Then a simple
calculation shows that $B^*A$ is
$$
\left(
\matrix
2 \Delta & 0& 0 & a_2 \Delta\\
0 & 2 \Delta & 0 & a_5  \Delta \\
0 & 0 & 2 \Delta & -b_1 \Delta
\endmatrix
\right)
$$

Assume now that $h$ is irreducible on the 3-holed sphere
bounded by $\alpha_2$ and $\alpha_5$, i.e., $\Delta \neq 0$.
Then we obtain a simpler system
of linear equations satisfied by $x_{12}, x_{15}, x_{34}$ and $x_1$.

$$
\aligned
 2 x_{12} + a_2 x_1 = p_{14} \\
 2 x_{15} + a_5 x_1 = p_{15} \\
 2 x_{34} - b_1 x_1 = p_{16}. 
\endaligned\tag 12
$$

By equations (1), (2), (3) and (12), we obtain
$$
\aligned
2 x_{21} +  a_2 x_1 = p_{17} \\
2x_{51} + a_5  x_1 = p_{18} \\
2x_{43} - b_1 x_1 = p_{19}.
\endaligned\tag 13
$$
On the other hand, there
are many different extensions of $\{\alpha_2, \alpha_5\}$
to a set forming a pentagon. For instance, $\{\alpha_1',
\alpha_2, \alpha_3', \alpha_4', \alpha_5\}$ =\{$ \alpha_1 \alpha_2,
\alpha_2, \alpha_3 \alpha_2, \alpha_4, \alpha_5\}$ is
such an extension. For this extension, the same equations (12) and
(13) hold. Here we have $x'_{21} = h(\alpha_{2} \alpha_1')
= h(\alpha_2 \alpha_1 \alpha_2) = h(\alpha_1) = x_1$,
$x_{51}' = h(\alpha_5 \alpha_1') = h(\alpha_5 \alpha_1 \alpha_2)
= h(\alpha_3\alpha_4) = x_{34}$ and $x_1' = h(\alpha_1')  = x_{12}$.
By (13) for the new pentagon, we obtain 

$$\aligned
2 x_1 + a_2 x_{12} = p_{20} \\
2x_{34}  + a_5 x_{12} = p_{21}.
\endaligned\tag 14
$$

Comparing (14) with (12), we obtain

$$(4 -a_2^2) x_1 = p_{22}, \tag 15$$

and 

$$(2b_1 - a_2 a_5) x_1 = p_{23}. \tag 16$$

Apply the involution $\tau$ to (15), we obtain 

$$(4 -a_5^2) x_1 = p_{24}. \tag 17$$ 

Since 
$a_2^2 + a_5^2 + b_1^2 - a_2 a_5 b_1 -4 \neq 0$, if the 
characteristic of the field $K$ is not 2, then
one of the coefficients $ 4-a_2^2$, $4-a_5^2$ and $2 b_1 -a_2a_5$
is not zero. Thus we can solve $x_1$ uniquely
from (15)-(17). If the characteristic of the field $K$ is 2, then
one of the coefficients $a_2$, $a_5$ or $b_1$ is not zero. 
Thus we can solve $x_1$ uniquely from (12).

Now take $h= f$ and $h = g$ respectively. The condition $f|_{S(X) \cup
S(Y)} = g|_{S(X) \cup S(Y)}$ shows that $x_1 = f(\alpha_1)$ and $x_1
= g(\alpha_1)$ are the solutions of the \it same \rm equations (1) -(17).
Since both $f$ and $g$ are irreducible on the 3-holed sphere bounded
by $\alpha_2$ and $\alpha_5$, we conclude that $f(\alpha_1) = g(\alpha_1)$.
Thus by corollary 3.7, $f= g$ follows in this case.

4.4.
In this section, we prove that if $f$ is a trace function which is
reducible on each level-0 subsurface, then $f$ is the character of
a reducible representation.

We choose a 3-holed sphere decomposition of $\Sigma_{0,5}$ by $\alpha_2$ and
$\alpha_5$ as follows. If there exists $\alpha \in S'(\Sigma_{0,5})$ so that
$f^2(\alpha) \neq 4$ then choose $\alpha_2$ to be one of these
elements. If otherwise that $f^2(\alpha) =4$ for all $\alpha \in S'(\Sigma_{0,5})$,
choose $\alpha_2$ and $\alpha_5$ to be any pair of disjoint elements.
We shall use the same notations introduced in \S4.3. Thus $X$ and $Y$
are level-1 subsurfaces bounded by $\alpha_2$ and $\alpha_5$
respectively so that $X \cap Y$ is a level-0 subsurface. By proposition
3.5, we find two representations $\rho_X$ and $\rho_Y$ of
$\pi_1(X)$ and $\pi_1(Y)$ respectively so that their characters are
the restrictions of $f$ to $S(X)$ and $S(Y)$. By the reducibility
criterion lemma 2.5, both $\rho_X$ and $\rho_Y$ are reducible.
Thus we may modify $\rho_X$ and $\rho_Y$ without changing their
characters so that both $\rho_X(\pi_1(X \cap Y))$ and $\rho_Y(\pi_1(X \cap Y))$
consist of diagonal matrices.  Now since both $\rho_X |_{\pi_1(X \cap Y)}$
and $\rho_Y|_{\pi_1(X \cap Y)}$ are diagonalizable and have the
same character, thus they are conjugate. We may assume after a
conjugation that $\rho_X|_{\pi_1(X \cap Y)} = \rho_Y|_{\pi_1( X \cap Y)}$. 
By the
same argment as in \S4.3, we construct a diagonalizable  representation $\rho$
of $\pi_1(\Sigma_{0,5})$ to $SL(2, K)$ extending both $\rho_X$ and
$\rho_Y$. Let $g$ be the character of $\rho$ defined on $S(\Sigma_{0,5})$. By the
construction $f|_{S(X)  \cup S(Y)}  = g|_{S(X) \cup S(Y)}$. The
goal is to show that $f= g$ under the reducible condition.

Since $f$ is reducible on all level-0 subsurfaces, by corollary 2.5,
$f$ is reducible on  all level-1 subsurfaces. In particular, by
corollary 3.5(d),  $f(\alpha \beta) = f(\beta \alpha)$ for all $\alpha \perp_0 \beta$. 

We now set up the same system of linear equations in $x_{i_1...i_k}$, 
$1 \leq k \leq 3$ as in \S4.3. Then equations (1) - (10) still hold. Due to
the reducibility, $x_{12} = x_{21}$, $x_{15} = x_{51}$ and
$x_{34} = x_{43}$. Thus equations (12)-(17) still hold. (Indeed,
equation (12) is a consequence of (1) and $x_{12}= x_{21}$.) 

Now if there is $\alpha \in S'(\Sigma_{0,5})$ so that
$f^2(\alpha) \neq 4$, then $4-f^2(\alpha_2) \neq 0$ by the choice
of $\alpha_2$. Thus we can solve $x_1$ uniquely from (15). In
particular, by the same argument as in \S4.3, we obtain $f(\alpha_1)
= g(\alpha_1)$.

In the remaining case, $f^2(\alpha) =4$ for all $\alpha \in S'(\Sigma_{0,5})$.
First we note that corollary 3.3(a) implies 
$f^2(\beta) = 4$ for all boundary
component $\beta$ in $S(\Sigma_{0,5})$. 
For each boundary component $\beta$, 
$f(\beta) = g(\beta)$ by the construction. We claim that
that $f(\alpha)
=g(\alpha)$ for all $\alpha \in S'(\Sigma_{0,5})$.  Indeed, for
each $\alpha \in S'(\Sigma_{0,5})$, there exists two boundary components
$\beta_1$ and $\beta_2$ so that $\beta_1, \beta_2$ and $\alpha$ bound
a level-0 subsurface. By the reducibility of $f$ and $g$ over
the level-0 subsurface and corollary 3.3(a), we conclude that
$f(\alpha) = f(\beta_1)f(\beta_2)/2 = g(\beta_1)g(\beta_2)/2
= g(\alpha)$. (Here we have used the convention that if the
characteristic of $K$ is 2, then $ab/2$ is meant to be $b$ when $a=2$).

4.5. As a consequence of theorem 4.1, we have,

{\bf Corollary.} \it
Let $f$ be a trace function defined on $S(\Sigma_{0,n})$.
Suppose $\Sigma_{0,n}$ is decomposed as a union  $X_1 \cup X_2 $ of two
incompressible subsurfaces $X_1$ and $X_2$
where $X_1 \cap X_2 \cong \Sigma_{0,3}$ is bounded by two elements in
$S'(\Sigma)$.
If $f|_{S(X_i)}$ is an $SL(2,K)$ character for $i=1,2$ and
either $f|_{S(X_1 \cap X_2)}$ is irreducible or $f$ is 
reducible on  all level-0 subsurfaces, then $f$ is an $SL(2,K)$
character. \rm

\it Proof. \rm Let $\rho_i$ be an $SL(2,K)$ representation
whose character is $f|_{S(X_i)}$, $i=1,2$.
If $f$ is reducible on all level-0 subsurfaces, then
by corollary 3.5(c), both $\rho_1$  and $\rho_2$ are reducible. In this
case, we may assume without changing the characters that both
$\rho_1$ and  $\rho_2$ are diagonalizable.
Now by the same
arguments as in \S4.3 and \S4.4, we  produce a representation
$\rho$ of $\pi_1(\Sigma_{0,n})$ so that its restriction to
$\pi_1(X_i)$ is  conjugate to $\rho_i$. Let $g$ be the character
of $\rho$ and $\beta_1$ and $\beta_2$ be two classes
in $S'(\Sigma_{0,n})$ which bound $X_1 \cap X_2$. 
Then $f$ and $g$ are identical on $S(X_1)
\cup S(X_2)$ =\{$ \alpha \in S(\Sigma_{0,n}) | $ $\alpha$ is
disjoint from either $\beta_1$ or $\beta_2$\}.
To show $f=g$, by corollary 3.7, it suffices to prove
$f(\alpha) = g(\alpha)$ for $\alpha \perp_0 \beta_i$, $i=1,2$.
Fix such a class $\alpha$. Let $\Sigma'$ be the incompressible
level-2 subsurface which contains $X_1 \cap X_2$ and $\alpha$.
Then  by the proof of theorem  4.1 for $\Sigma'$ with
respect to the decomposition $\Sigma' = (\Sigma' \cap X_1) \cup
(\Sigma' \cap X_2)$, 
we have
$f|_{S(\Sigma')} = g|_{S(\Sigma')}$. In particular, $f(\alpha)
= g(\alpha)$. 
$\square$

\bigskip

{\bf \S5. Exceptional trace functions on Planar Surfaces}

5.1. Recall that a trace function which is not the character of
any representation is called \it exceptional. \rm There are no
exceptional trace functions on level-0 and level-1 surfaces. However,
there exist finitely many exceptional trace functions on $\Sigma_{0,n}$
for any $n \geq 5$. The main result of the section is to identify
all exceptional trace functions.

{\bf Theorem.} \it Suppose $f$: $S(\Sigma_{0,n}) \to K$, $n \geq 5$,
is an exceptional trace function. Then the characteristic of $K$
is not 2 and $f$ satisfies,

(a) $f(S(\Sigma_{0,n})) =\{2, -2\}$ and, 

(b) there exists an exceptional level-2 subsurface in $\Sigma_{0,n}$.

\rm 

The proof of the theorem is by induction on $n$. In \S5.2 we prove it
for $n=5$ and in \S5.3, we prove it for all $n \geq 6$. 

We shall use the following notations. If $\alpha_1, ..., \alpha_m $
are disjoint classes in  $S'(\Sigma_{g,n})$
so that they decompose the surface into subsurfaces $\Sigma_{g_i, n_i}$
and $(g_1, n_1) \neq (g_i, n_i)$ for $i \geq 2$, then
we use $\Sigma_{g_1, n_1}(\alpha_1, ..., \alpha_m)$ to
denote the subsurface $\Sigma_{g_1, n_1}$.
A  class  $\alpha \in S'(\Sigma)$
is a \it boundary class \rm if it bounds a level-0 subsurface.

5.2. We prove  a slightly stronger version of theorem 5.1 for $n=5$ in this 
section. 

Let $b_1, ..., b_5$ be the boundary components of $\Sigma_{0,5}$.

{\bf Proposition.} \it If $f: S(\Sigma_{0,5}) \to K$ is an exceptional
trace function, then

(a) the characteristic of $K$ is not 2, 

(b) $f(S(\Sigma_{0,5})) =\{2, -2\}$ and $\Pi_{i=1}^5 f(b_i) = 32$, and,

(c) a level-0
subsurface is irreducible if and only if it is of the form $\Sigma_{0,3}
(\alpha)$ for some $\alpha \in S'(\Sigma_{0,5})$. \rm

\it Proof. \rm By theorem 4.1, we see that $\Sigma_{0,3}(\alpha, \beta)$
is reducible for all disjoint $\alpha, \beta \in S'(\Sigma_{0,5})$ 
and there exists
one irreducible $\Sigma_{0,3}(\gamma)$.

{\bf Lemma.} \it Suppose $\alpha_1, \alpha_4$ are two disjoint elements
in $S'(\Sigma_{0,5})$ so that $\Sigma_{0,3}(\alpha_1)$ is irreducible.
Then  $f^2(\alpha_1) = f^2(\alpha_4) = f^2(b_i) = 4$ and $\Sigma_{0,3}(\alpha_4)$
is again irreducible. Furthermore, the characteristic of $K$ is not 2.
\rm

\it Proof. \rm Since $\Sigma_{0,3}(\alpha_1)$ $\subset \Sigma_{0,4}(\alpha_4) $
and $\Sigma_{0,3}(\alpha_1)$  is irreducible, we see that
$\Sigma_{0,4}(\alpha_4)$ is  irreducible. By corollary 3.5(a) applied to
$\alpha_1$ in $\Sigma_{0,4}(\alpha_4)$ and by the assumption that
$\Sigma_{0,3}(\alpha, \alpha_4)$ is reducible, there exists $\alpha_2
\perp_0 \alpha_1$ in $\Sigma_{0,4}(\alpha_4)$ so that $\Sigma_{0,3}(\alpha_2)$
is irreducible. Extend $\{\alpha_1, \alpha_2, \alpha_4\}$ to a pentagon set
$\{\alpha_1, \alpha_2, \alpha_3, \alpha_4, \alpha_5\}$ so that
$\alpha_i \perp_0 \alpha_{i+1}$ where indices are  counted mod 5. For
$i=3,4,5$, $\Sigma_{0,4}(\alpha_i)$ is irreducible since it contains
one of $\Sigma_{0,3}(\alpha_j)$, $j=1$ or $2$. By corollary 3.5(b) applied
to $\alpha_i$ in $\Sigma_{0,4}(\alpha_i)$ and by the 
assumption that $\Sigma_{0,3}(\alpha_i, \beta)$ is reducible, it follows that $f^2(\alpha_i) =4$
for $i=3,4,5$. Let the boundary components of $\Sigma_{0,5}$ be
so labelled that $(\alpha_{i-1}, b_i, \alpha_{i+1})$ bounds
a 3-holed sphere. For each $i=1,2,..., 5,$ $\Sigma_{0,3}(\alpha_{i-1},
\alpha_{i+1})$ is reducible and one of $f^2(\alpha_{i-1})$ or
$f^2(\alpha_{i+1})$ is $4$. By corollary 3.3(a), it follows that
$f(b_i) = f(\alpha_{i-1})f(\alpha_{i+1})/2$ (here  $ab/2$ is meant to
be $b$ if the characteristic of $K$ is $2$ and $a = 2$). In 
particular $f^2(b_4) = 4$. 
The values of $f$ on
$\partial \Sigma_{0,3}(\alpha_1)$ are \{$f(\alpha_1), f(b_3), f(b_4)
$\} 
$= \{f(\alpha_1), f(\alpha_2)f(\alpha_4)/2,$$ 
f(\alpha_3)f(\alpha_5)/2$\} whose  multiplication is $\frac{1}{4}
\Pi_{i=1}^5 f(\alpha_i)$. For $i=3,4,5$, the values of $f$ 
on $\partial \Sigma_{0,3}(\alpha_i)$
are \{$f(\alpha_i)$, $f(b_{i+2})$, $f(b_{i+3})$\}
 $= \{f(\alpha_i), $$f(\alpha_{i-1})f(\alpha_{i+1})/2,$$
f(\alpha_{i+2})f(\alpha_{i+3})/2$\} whose multiplication is again 
$\frac{1}{4} \Pi_{i=1}^5 f(\alpha_i)$.
Since $\Sigma_{0,3}(\alpha_1)$ is 
irreducible and $f^2(b_4) = f^2(\alpha_i)
= 4$ for $i=3,4,5$, by corollary 3.3(a),
it  follows that $\Sigma_{0,3}(\alpha_i)$
is irreducible. 

Now by the same argument above applied to $\{\alpha_3, \alpha_4\}$
instead of $\{\alpha_1, \alpha_2\}$, we conclude that $f^2(\alpha_1)
= f^2(\alpha_2) =4$. Thus $f^2(b_i) =4$ for all $i$. In particular,
this implies that the characteristic of $K$ is not $2$. Indeed,
if otherwise, then all $f(\alpha_i) = f(b_i) = 0$. Thus $\Sigma_{0,3}(\alpha_1)$ would be reducible.
$\square$

To finish the proof of the proposition, 
take $\alpha_1 \in S'(\Sigma_{0,5})$ so that
$\Sigma_{0,3}(\alpha_1)$ is irreducible. Given any 
$\alpha \in S'(\Sigma_{0,5})$,
by a result of  Harvey [Hav] (see also [Har]), there  exists a sequence of
elements $\beta_1 = \alpha_1, \beta_2, ..., \beta_m = \alpha$ in
$S'(\Sigma_{0,5})$ so that $\beta_i \cap \beta_{i+1} = \emptyset$.
By the lemma applied to the sequence, we conclude that $f^2(\alpha) =4$
and $\Sigma_{0,3}(\alpha)$ is irreducible.

Suppose $\alpha \in S'(\Sigma_{0,5})$ so that $b_i$ and $b_j$
are in the boundary of $\Sigma_{0,3}(\alpha)$. Then by corollary
3.3(a), $f(\alpha) = -\frac{1}{2}f(b_i)f(b_j)$. This shows that
$f$ is determined by $f|_{\partial \Sigma_{0,5}} $. Furthermore,
take two disjoint $\alpha,\alpha'$ in $S'(\Sigma_{0,5})$ so that
$b_1 \subset \partial \Sigma_{0,3}(\alpha, \alpha')$.
Then due to the reducibility and corollary 3.3(a), 
$f(b_1)f(\alpha)f(\alpha') =8$. But
also $f(\alpha) f(\alpha') = \frac{1}{4}f(b_2) f(b_3)f(b_4) f(b_5)$.
Thus $\Pi_{i=1}^5 f(b_i) = 32$.
$\square$

5.3. We use the induction on $n$ to prove theorem 5.1. Assume that
$n \geq 6$. The proof consists of several steps.

Let $f$ be an exceptional trace function on $S(\Sigma_{0,n})$.

Claim 1. \it There exists an exceptional subsurface $\Sigma_{0,n-1}(\beta)$
in $\Sigma_{0,n}$. \rm

By corollary 4.5, there exists an irreducible level-0 subsurface
$\Sigma'$ in $\Sigma_{0,n}$. Now $\Sigma'$ is either $\Sigma_{0,3}(\alpha_1,
\alpha_2)$, or $\Sigma_{0,3}(\alpha_1)$, or $\Sigma_{0,3}(\alpha_1, \alpha_2,
\alpha_3)$ for disjoint classes $\alpha_1, .., \alpha_3$
in $S'(\Sigma_{0,n})$. If $
\Sigma' = \Sigma_{0,3}(\alpha_1, \alpha_2)$, by corollary 4.5, one of
the subsurface bounded by $\alpha_1$ or $\alpha_2$ is exceptional. Take
any $\Sigma_{0, n-1}(\beta)$ which contains this exceptional subsurface,
then the  claim follows.

If $\Sigma' = \Sigma_{0,3}(\alpha_1)$, we use the following lemma.

{\bf Lemma.} \it If $\Sigma_{0,3}(\alpha)$ is irreducible
and $\beta_1, \beta_2$ are two disjoint classes in $S'(\Sigma_{0,n})$
so that they bound a level-1 subsurface which contains $\Sigma_{0,3}(\alpha)$,
then one of the  subsurface bounded by $\beta_i$ is exceptional. \rm

\midspace{0.1cm}
\centerline{\epsfbox{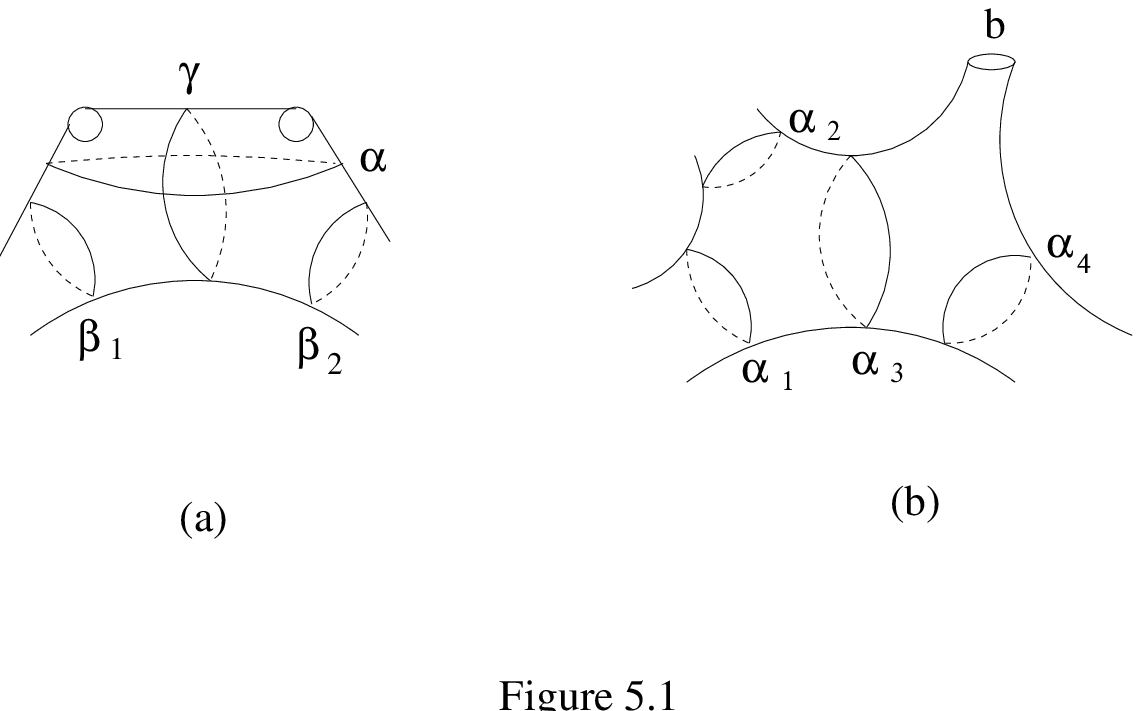}}

Indeed, by corollary 3.5(a), there exists $\gamma \perp_0 \alpha$
in $\Sigma_{0,4}(\beta_1, \beta_2)$ so that one of $\Sigma_{0,3}(\gamma,
\beta_i)$, say $\Sigma_{0,3}(\gamma, \beta_1)$, is
 irreducible. By corollary 4.5,
one of the subsurface bounded by $\gamma$ or $\beta_1$ is exceptional.
But each subsurface bounded by $\gamma$ or $\beta_1$ is contained
in a subsurface bounded by $\beta_1$ or $\beta_2$. Thus the lemma follows.

If $\Sigma' = \Sigma_{0,3}(\alpha_1, \alpha_2, \alpha_3)$, let
$\Sigma''$ be the  irreducible subsurface bounded by $\alpha_1, \alpha_2$
which contains $\Sigma'$. If $\Sigma''$ is exceptional, then claim 1 follows.
If the $\Sigma''$ is a 4-holed sphere, then by corollary 3.5(a)
applied to $\alpha_3$ in $\Sigma''$, we find an  irreducible
$\Sigma_{0,3}(\alpha_i, \alpha')$ where $i=1$ or $2$.  The claim
follows by the previous argument.
In the remaining case,  the level 
of $\Sigma''$ is at least 2 and  $f|_{S(\Sigma'')}$ 
is the character of a representation
$\rho$ of $\pi_1(\Sigma'')$. Since $\rho$ is irreducible, there
exists a boundary component $b \subset \partial \Sigma_{0,n} \cap
\partial \Sigma''$ so that $\rho(b) \neq \pm id$. Consider a level-1
subsurface of $\Sigma''$ of the form $\Sigma_{0,4}(\alpha_1, \alpha_2,
\alpha_4)$ so that it contains $\Sigma_{0,3}(\alpha_1, \alpha_2,
\alpha_3)$ and contains $b$ as a boundary component. By
corollary 3.5(b) applied to this level-1 subsurface, we find 
an irreducible  $\Sigma_{0,3}(\alpha, \alpha')$ having $b$
as a boundary component.
Thus the claim follows.

Now by the induction hypothesis applied to $\Sigma_{0, n-1}(\beta)$,
it follows that $\Sigma_{0,n}$ contains an exceptional level-2 subsurface,
i.e., part (b) of the theorem 5.1 follows.

Claim 2. \it Let $b_1$ and $b_2$ be the boundary components of $\Sigma_{0,n}$
which are not in the exceptional $\Sigma_{0,n-1}(\beta)$. Then $f^2(b_i)= 4$
for $i=1,2$. \rm

Indeed, if $n \geq 7$, for each $b_i$, there exists an $(n-1)$-holed
sphere containing both $b_i$ and the exceptional level-2 subsurface
in $\Sigma_{0,n-1}(\beta)$. Thus by the induction  hypothesis, we
conclude that $f(b_i)^2 = 4$ for $i=1,2$.

If $n=6$, we pick a boundary class $\beta' \in S'(\Sigma_{0,6})$ so that
$\Sigma_{0,3}(\beta')$  is in
the exceptional  $\Sigma_{0,5}(\beta)$ (see fig. 5.2). 
Let $\beta_1$ and $\beta_2$ be two disjoint boundary classes
so that $\beta_i \cap \beta' = \emptyset$ and $\beta_i \perp_0
\beta$. Since $\Sigma_{0,5}(\beta)$ is exceptional, $\Sigma_{0,3}(\beta')$
is irreducible. By the above lemma applied to $\alpha = \beta'$,
we conclude that one of $\Sigma_{0,5}(\beta_i)$, say $\Sigma_{0,5}(\beta_1)$,
is exceptional. Since $b_1 \subset \Sigma_{0,5}(\beta_1)$, it follows
that $f^2(b_1) = 4$. Next we assert that $\Sigma_{0,3}(\beta)$ 
is reducible. Assuming
the  assertion, by corollary 3.3(a), we conclude that $f(b_2)
= f(b_1) f(\beta)/2 = \pm 2$.
To see that $\Sigma_{0,3}(\beta)$ is reducible, we construct 
two disjoint boundary classes $\gamma_1, \gamma_2$ in
$\Sigma_{0,6}$ which are in $\Sigma_{0,5}(\beta)$. If
$\Sigma_{0,3}(\beta)$ were irreducible, then by the above lemma 
applied to $\{\beta, \gamma_1, \gamma_2\}$, we conclude that
one of $\Sigma_{0,5}(\gamma_i)$, say $\Sigma_{0,5}(\gamma_1)$,
is exceptional. Let $\gamma$ be a class disjoint from 
$\beta$ and $\gamma_1$ and $\gamma \perp_0 \gamma_2$. Then
$\Sigma_{0,3}(\beta, \gamma)$ is irreducible since it
is in the exceptional $\Sigma_{0,5}(\beta)$. But $\Sigma_{0,3}(\beta,
\gamma)$ is also reducible since it is in the exceptional $\Sigma_{0,5}
(\gamma_1)$. This is a contradiction.

\midspace{0.1cm}
\centerline{\epsfbox{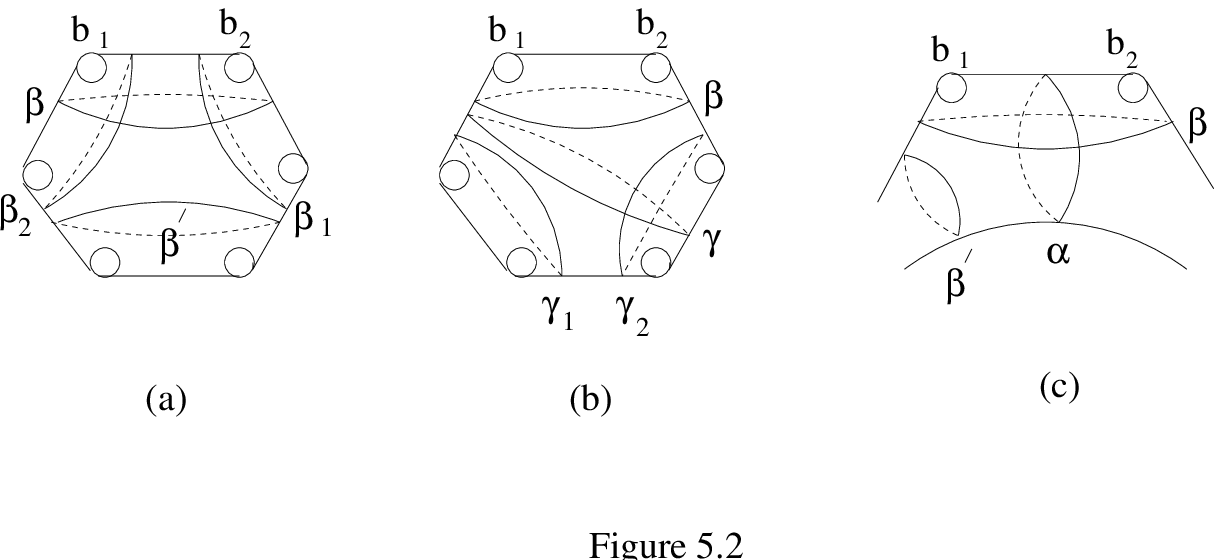}}

Claim 3. \it If $\alpha \perp_0 \beta$ where $\Sigma_{0,n-1}(\beta)$
is exceptional, then $f(\alpha) = \pm 2.$ \rm

Choose a class $\beta' \in S'(\Sigma_{0,n})$ 
disjoint from $\beta$ so that
$\beta'$, $\alpha$ and one of $b_1$ or $b_2$, say $b_1$,
bound a 3-holed
sphere $\Sigma_{0,3}(\beta', \alpha)$ (see fig. 5.2(c)). 
By the induction hypothesis
applied to $\Sigma_{0,n-1}(\beta)$, we have $f(\beta') = \pm 2$.
If $\Sigma_{0,3}(\alpha, \beta')$ is reducible, then
by corollary 3.3(a), we obtain $f(\alpha) = f(\beta') f(b_1)/2
=\pm 2$.
If $\Sigma_{0,3}(\alpha, \beta')$ is irreducible, then
by corollary 4.5, one of the subsurface $X$ bounded by $\beta'$ or
$\alpha$ which contains $\Sigma_{0,3}(\alpha, \beta)$ is
exceptional. Thus by the induction hypothesis applied
to a subsurface containing $X$ and $\alpha$, $f(\alpha) = \pm 2$.

Claim 4. \it For all $\alpha \in S'(\Sigma_{0,n})$, $f(\alpha) = \pm 2$.
\rm

We use the induction on $I(\alpha, \beta)$ to prove the claim.
By claim 3, the result holds for
$I(\alpha, \beta) \leq 2$. If $I(\alpha, \beta) \geq 4$, by lemma 3.7,
we can express $\alpha = \alpha' \alpha''$ so that $\alpha'
\perp_0 \alpha''$, $\partial (\alpha', \alpha'') = \beta_1 \cup ..
. \cup \beta_4$ satisfy $I(\alpha', \beta), I(\alpha'', \beta), 
$ $ I(\alpha'' \alpha', \beta), 
I(\beta_i, \beta) <  I(\alpha, \beta)$. Thus by the induction
hypothesis the values of $f$ on the seven elements $\{
\alpha', \alpha'', \alpha'' \alpha', \beta_1, ..., \beta_4\}$ 
are in $\{2, -2\}$. Now the following lemma implies that $f(\alpha)
= \pm 2.$

5.4.
{\bf Lemma.} \it Let $\partial \Sigma_{0,4} = b_1 \cup b_2 \cup b_3
\cup b_4$ and $\alpha_1$, $\alpha_2$, and $\alpha_3$ be three classes
forming a triangle in $S'(\Sigma_{0,4})$. If $f: S(\Sigma_{0,4}) \to
K$ is a character so that its values on the 7-element set $\{\alpha_i,
b_j\}$ are $\{2, -2\}$. Then $f(S(\Sigma_{0,4})) \subset \{2, -2\}$.
Furthermore,

(a) $2\Pi_{i=1}^3 f(\alpha_i) = \Pi_{j=1}^4 f(b_j)$,

(b) if the characteristic of $K$ is not $2$, then $f$ is reducible
if and only if $\Pi_{j=1}^4 f(b_j) = 16$.

(c) If $g: \{\alpha_i, b_j\} \to \{2, -2\}$ satisfies $2\Pi_{i=1}^3 
g(\alpha_i ) = \Pi_{j=1}^4 g(b_j) = -16$, then $g$ can be extended to
an $SL(2, K)$ character on $S(\Sigma_{0,4})$. \rm

\it Proof. \rm Fix an orientation on each $b_i$ and consider
it as an element in the fundamental group. Let $\rho$ be a
representation whose character is $f$. Changing $\rho$ to $\rho'$
  by $\rho'(b_i) = \pm \rho(b_i)$, 
$i=1,2,3$, will not effect the conclusion of the lemma. Thus
 we may assume that $f(b_1) = f(b_2) = f(b_3) = 2$. Now if
$f(b_4) = 2$, then proposition 3.5(a) shows that
$8\Sigma_{i=1}^3 f(\alpha_i) - \Pi_{i=1}^3 f(\alpha_i) = 40$.
Since $f(\alpha_i) = \pm 2$, the only solution of the equation
is $f(\alpha_i) = 2$. By proposition 3.5(b), this implies that
$f(\alpha) = 2$ for all $\alpha$. If $f(b_4) = -2$, then proposition
3.5(a) says that $\Pi_{i=1}^3 f(\alpha_i) = -8$.
Furthermore, proposition 3.5(b) implies that $f(\alpha_3') + f(\alpha_3)
= -f(\alpha_1) f(\alpha_2)$ where $(\alpha_1, \alpha_2, \alpha_3;
\alpha_3')$ forms a quadrilateral. But $f(\alpha_3) = - \frac{1}{2}f(\alpha_1)
f(\alpha_2)$. Thus $f(\alpha_3') = f(\alpha_3) = \pm 2$. By the
modular configuration, this implies that $f(S(\Sigma_{0,4})) \subset
\{2, -2\}$.  The last argument also shows that for any assignments
of $\pm 2$ to $\alpha_i$'s so that their product is $-8$, there
exists an extension of the assignment to a character. Thus part (c)
follows. 
$\square$

The following figure 5.3 illustrates the set of all possible assignments
of $\pm 2$ to the 7-element set $\{b_i ,\alpha_j\}$ in the lemma.

\midspace{0.1cm}
\centerline{\epsfbox{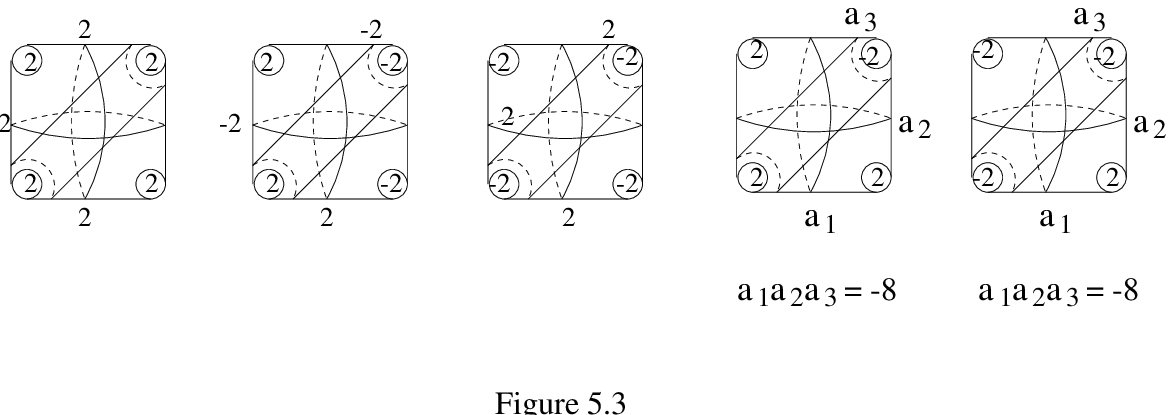}}

As a consequence of the  theorem 5.1(b)  and the fact that trace functions
are determined by their values on a finite subset of $S(\Sigma)$ (\S3.9),
we
see that there are only finitely many exceptional trace functions
on each planar surface.

5.5. 
The goal of this section is to prove that exceptional trace functions
exists on each planar surface of level at least 2.

There are sixteen exceptional trace functions $f$ on $\Sigma_{0,5}$ which we
describe as follows. Let $\partial \Sigma_{0,5} = b_1 \cup ... \cup b_5$.
Suppose $f: \{b_1, ..., b_5\} \to \{2, -2\}$ satisfies 
$\Pi_{i=1}^5 f(b_i) = 32$. We extend $f$ to $f: S(\Sigma_{0,5})
\to \{2, -2\}$ as follows. Given any
boundary class $\alpha$ so that the level-0
subsurface $\Sigma_{0,3}(\alpha)$ contains
$b_i, b_j$, then $f(\alpha) = -\frac{1}{2}f(b_i) f(b_j)$.
Ones checks easily using lemma 5.4 that $f$ is a trace function
on $S(\Sigma_{0,5})$. Furthermore, by the construction
$f$ is reducible on each $\Sigma_{0,3}(\alpha, \beta)$ and
irreducible on each $\Sigma_{0,3}(\alpha)$.
There is no representation whose character is the
trace function $f$. 
Indeed, if $\chi_{\rho} = f$  for a representation $\rho$
and $\alpha \in S'(\Sigma_{0,5})$, then
by corollary 3.5(b) applied to $\Sigma_{0,4}(\alpha)$, 
we have $\rho(\alpha) = \pm id$.
But this implies that $\rho$ is reducible on $\Sigma_{0,3}(\alpha)$
which contradicts the assumption.
Thus these are the set of all 16 exceptional trace functions on
$\Sigma_{0,5}$.

If $n \geq 6$, we construct an exceptional trace function on $\Sigma_{0,n}$
as follows. Let $b_1, ..., b_n$ be the boundary components of $\Sigma_{0,n}$.
Define $f: S(\Sigma_{0,n}) \to \{2, -2\}$ as follows.
Let $f(b_i) = 2$ for all $i$.
For $\alpha \in S'(\Sigma_{0,n})$, we define $f(\alpha)$ as follows.
Suppose $\alpha$ decomposes $\Sigma_{0,n}$ into two subsurfaces
$X_1$ and $X_2$. Let $s_i$ be the number of components of $\{b_1, ..., b_5\}$
which are in $X_i$. Define $f(\alpha)$ to be $-2$ if $(s_1, s_2) = (2,3)$
and to be $2$ otherwise.  By a simple calculation using lemma 5.4, one
shows that $f$ is a trace function and there exists an exceptional
level-2 subsurface. Thus $f$ is an exceptional trace function.

\S6. {\bf The Characterization Theorem for 1-holed Torus}

We show that each trace function on $S(\Sigma_{1,2})$ is a character in this
section.

6.1. \it The pentagon relation \rm

A five-element set \{$\alpha_1, ..., \alpha_5\}$ in $S'(\Sigma_{1,2})$ is said
to form a \it pentagon \rm if $\alpha_i \cap \alpha_{i+2} =
\emptyset$ for all $i$ mod 5. It is shown in [Lu3] that the
set is unique up to homeomorphism of the surface and
that exactly two adjacent elements $\alpha_i$ and $\alpha_{i+1}$, say
$\alpha_3$ and $\alpha_4$, are separating classes and
$\alpha_3 \perp_0 \alpha_2 \perp \alpha_1 \perp \alpha_5 \perp_0
\alpha_4$. See figure 6.1.

If \{$\alpha_1,..., \alpha_5\}$ forms a pentagon with $I(\alpha_3, \alpha_4)
= 4$, then we have,

(a)  $ (\alpha_i \alpha_j) \alpha_k = \alpha_i(\alpha_j \alpha_k) $
where the indices $i,j,k$ are pairwise distinct,

(b) $\alpha_i \alpha_j \alpha_k = \alpha_j \alpha_i \alpha_k$ if
the indices are pairwise distinct and $\alpha_i \cap \alpha_j = \emptyset$,

(c) $\alpha_2 (\alpha_2 \alpha_1 \alpha_5) = \alpha_4 \alpha_5 \alpha_1$,
 $(\alpha_2 \alpha_1 \alpha_5 ) \alpha_2 =  \alpha_1 \alpha_5 $,
and $\alpha_2 \alpha_1 \alpha_5 \cap \alpha_1 = \emptyset$,


(d) $\alpha_1 \alpha_2 \cap \alpha_1 \alpha_5 = \emptyset$ and
$\alpha_2 \alpha_1 \cap \alpha_5 \alpha_1 =\emptyset$, 

(e) $\alpha_1 \alpha_2 \alpha_3 \cap \alpha_1 \alpha_5 = \emptyset$ and
$\alpha_4 \alpha_5 \alpha_5 \cap \alpha_1 \alpha_5 = \emptyset$.

See [Lu2] or [Lu3] for a verification. One can also verify (a)-(e)
directly. For instance $\alpha_1 \alpha_5$ and $\alpha_1 \alpha_2$
are obtained by applying the positive Dehn twist along $\alpha_1$
to $\alpha_5$ and $\alpha_2$. Thus (d) holds.

By property (d), if
\{$\alpha_1, \alpha_2 , \alpha_3, \alpha_4, \alpha_5\}$
forms a pentagon, then 
\{$\alpha_1, \alpha_2 \alpha_1, \alpha_3, \alpha_4, \alpha_5\alpha_1\}$
is also a pentagon. 

6.2. In this section, we prove the following,

{\bf Proposition.} \it If  $f$ is a trace function on
$S(\Sigma_{1,2})$ and there exist two disjoint elements $\alpha_1, \alpha_4$
in $S'(\Sigma_{1,2})$ with $\alpha_4$ separating so that $f^2(\alpha_1)
\neq f(\alpha_4) + 2$, then $f$ is a character.
\rm

\it Proof. \rm
Let $X$ and $Y$ be the level-1 subsurfaces bounded by $\alpha_4$
and $\alpha_1$ respectively. Then $X \cap Y$ is a level-0 subsurface
bounded by $\alpha_4, \alpha_1, \alpha_1$. Since $f^2(\alpha_1)
+ f^2(\alpha_1) + f^2(\alpha_4) - f(\alpha_1) f(\alpha_1) f(\alpha_4)
-4 = (f(\alpha_4) -2) (f(\alpha_4) + 2 - f^2(\alpha_1))$, $f$ is reducible
on $X \cap Y$ if and only if $f(\alpha_4) = 2$.
We now construct a representation $\rho$ of $\pi_1(\Sigma_{1,2})$ so that
the  restrictions of the character of $\rho$ to 
$X$ and $Y$ are the same as $f|_{S(X)}$ and $f|_{S(Y)}$ as follows.
If $f(\alpha_4) \neq 2$, due to the irreducibility of $f$ on $X \cap Y$,
the construction is the same as in \S4.3. If $f(\alpha_4) = 2$, then
$f(\alpha_2) \neq \pm 2$ by the assumption. 
By lemma 2.6, there are  exactly two conjugation classes
of $SL(2,K)$ representations of $\pi_1(X)$ (respectively $\pi_1(X \cap Y)$)
whose characters are $f|_{S(X)}$ (resp. $f|_{S(X \cap Y)}$).
Furthermore, due to $f(\alpha_2) \neq \pm 2$,
 the restriction of the non-diagonalizable
representation of $\pi_1(X)$  to  $\pi_1(X \cap Y)$ is still
non-diagonalizable. Now  take a representation $\rho_Y$ of $\pi_1(Y)$ whose
character is $f|_{S(Y)}$. Then there exists a representation $\rho_X$ of
$\pi_1(X)$ whose  character is  $f|_{S(X)}$ so that $\rho_X|_{\pi_1(X \cap Y)}
= \rho_Y|_{\pi_1( X \cap Y)}$. 
Let $\rho$ be the representation  of $\pi_1(\Sigma)$ whose
restrictions to $\pi_1(X)$ and $\pi_1(Y)$ are $\rho_X$ and $\rho_Y$
and let $g$ be its character. We have $f(\alpha) 
= g(\alpha)$ for $\alpha \in S(X) \cup S(Y)$.
The goal is to show that $f=g$. By corollary 3.7, it suffices to
prove that for each class $\alpha_5 \perp \alpha_1$ and $\alpha_5
\perp_0 \alpha_4$, $f(\alpha_5) = g(\alpha_5)$.

Extend \{$\alpha_1, \alpha_4, \alpha_5\}$ to a 5-element set
$\{\alpha_1, ..., \alpha_5\}$ forming a pentagon. 
The proof of $f(\alpha_5) = g(\alpha_5)$ follows the same
strategy as in \S4.3 by introducing a system of linear equations.

We shall use the same notations as in \S4.3. Let $h$ be a 
trace function on $S(\Sigma_{1,2})$ so that $h^2(\alpha_1) \neq h(\alpha_4) +2$.
Given a set of indices $i_1, ..., i_k$, $1 \leq k \leq 3$, let
$x_{i_1...i_k} = h(\alpha_{i_1} ... \alpha_{i_k})$
if $\alpha_{i_1}...\alpha_{i_k}$ is not in $S(X) \cup S(Y)$ and
$a_{i_1...i_k} = h(\alpha_{i_1} ... \alpha_{i_k})$ if
$\alpha_{i_1}...\alpha_{i_k} \in S(X) \cup S(Y)$.
Let $\beta_1$ and $\beta_2$ be the boundary components of
$\Sigma_{1,2}$ and $b_i = h(\beta_i)$.

Using propositions 3.4 and 3.5, we now derive a system of linear
equations in $x_{i_1...i_k}$ and show that the system of equation
has a unique solution.

Since $\alpha_1 \perp \alpha_5$, by proposition 3.4(b), we obtain

$$h(\alpha_1 \alpha_5) + h(\alpha_5 \alpha_1) = h(\alpha_1) h(\alpha_5).$$

This is the same as,

$$ x_{15} + x_{51} = a_1 x_5. \tag 1$$

Since $\alpha_4 \perp_0 \alpha_5$ so that $\partial(\alpha_4, \alpha_5)
=\{\alpha_2, \alpha_2, \beta_1, \beta_2\}$, by proposition 3.5, we
obtain,

$$x_{45} + x_{54} = -a_4 x_5 + p_1.  \tag 2$$
where $p_1$ and the $p_i$'s below are polynomials with integer
coefficients in $h(\alpha)'s$ where $\alpha \in S(X) \cup S(Y)$.

Since $\alpha_2 \perp \alpha_2 \alpha_1 \alpha_5$ so that
$\alpha_2 (\alpha_2 \alpha_1 \alpha_5) = \alpha_4 \alpha_5 \alpha_1$
and $(\alpha_2 \alpha_1 \alpha_5) \alpha_2 = \alpha_1 \alpha_5$ and
$\alpha_2 \alpha_1 \alpha_5 \in S(Y)$,
 we obtain,

$$x_{451} + x_{15} = p_2. \tag 3$$

Let $\tau$ be the orientation reversing involution of $\Sigma_{1,2}$ fixing each
$\alpha_i$'s (see fig. 6.1). Then $\tau(\alpha \beta) = \tau( \beta) \tau(
 \alpha)$ for all $\alpha \perp \beta$ or $\alpha \perp_0 \beta$.

\midspace{0.1cm}
\centerline{\epsfbox{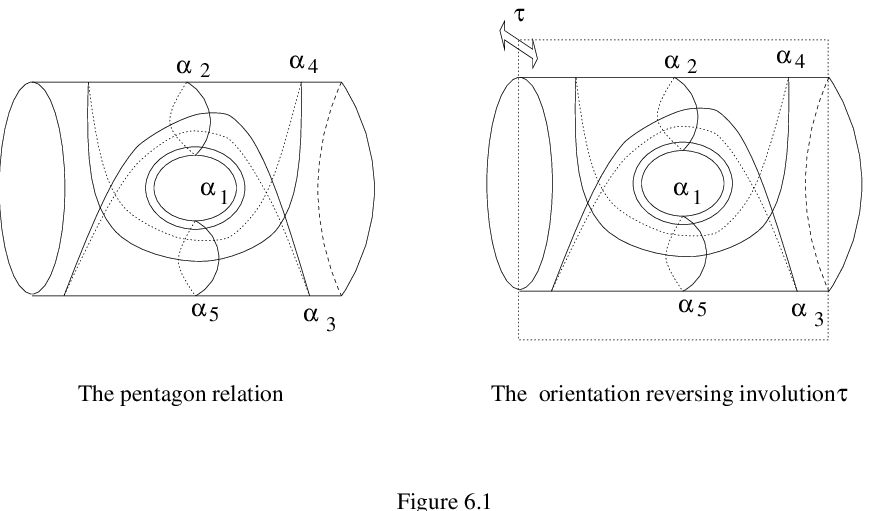}}

Apply $\tau$ to (3), we obtain
$$ x_{154} + x_{51} = p_3. \tag 4$$

Since $\alpha_4 \perp_0 \alpha_1 \alpha_5$ so that $\partial( \alpha_4, \alpha_1
\alpha_5) = \alpha_1 \alpha_2 \in S(X)$,  we obtain,

$$x_{415} + x_{154} + a_4 x_{15} = p_4. \tag 5$$

Since $\alpha_4 \cap \alpha_1 = \emptyset$,

$$ x_{415} = x_{145}. \tag 6$$

Finally, since $\alpha_1 \perp \alpha_4 \alpha_5$, by (6) we obtain,

$$x_{145} + x_{451} - a_1 x_{45} = 0. \tag 7$$

Subtracting (7) by (5) and using (6), we obtain

$$ x_{451} - x_{154} - a_1x_{45} - a_4x_{15} = p_5. \tag 8$$

The subtraction (4) by (3) gives,

$$x_{154} -x_{451} + x_{51} -x_{15} = p_6. \tag 9$$

The sum of (8) and (9) gives

$$x_{51} - x_{15} - a_1 x_{45} - a_4 x_{15} = p_7. \tag 10 $$

Using (1), we simplify (10) and obtain

$$a_1 x_{45} + (2 + a_4)x_{15} - a_1 x_5 = p_8. \tag 11$$

Since (11) holds for any 5-element set \{$\alpha_1', ..., \alpha_4', \alpha_5'\}$
forming a pentagon so that $\alpha_1' = \alpha_1$ and $\alpha_4'= \alpha_4$,
 it holds for the set
\{$\alpha_1, \alpha_2 \alpha_1, \alpha_3, \alpha_4, \alpha_5\alpha_1\}$
=\{$\alpha_1', \alpha_2', \alpha_3', \alpha_4', \alpha_5'\}$. Now
$ h(\alpha_1' \alpha_5') = h(\alpha_1 \alpha_5 \alpha_1) = h(\alpha_5)
= x_5$, $h(\alpha_5') = h(\alpha_5 \alpha_1) = x_{51} = a_1 x_5 -x_{15}$
and $h(\alpha_4' \alpha_5') = h(\alpha_4 \alpha_5 \alpha_1) =x_{451}
= -x_{15} + p_2$. Thus equation (11) for this new pentagon set
gives,

$$a_1(-x_{15} + p_2) + (2+a_4) x_5 -a_1(a_1 x_5 - x_{15}) = p_9$$

which is 

$$(2 + a_4 -a_1^2)x_5 = p_9. $$

Thus $x_5$ can be solved uniquely.  Now take $h=f$ and $h = g$. We
see that $f(\alpha_5) = g(\alpha_5)$. By corollary 3.7, it follows 
that $f=g$.

6.3. Suppose now that $f$  is a trace function so that
$f^2(\alpha_1) = f(\alpha_4) + 2$ for all separating $\alpha_4$ and
non-separating $\alpha_1$ with $\alpha_1 \cap \alpha_4 = \emptyset$.

We begin with the following lemma.

{\bf Lemma.} \it Let $\Sigma =\Sigma_{g,n}$ be a surface of level at
least 2 so that $g \geq 1$ and $f$ is a $K$-valued trace function
on $S(\Sigma)$. Let $P(\Sigma) =\{ (\alpha, \beta)
\in S(\Sigma) |$$\alpha$ bounds a $\Sigma_{1,1}$ in  $\Sigma$
and $\beta$ is a non-separating class lying in $\Sigma_{1,1}$\}.

(a) If for all $(\alpha, \beta) \in$ $P(\Sigma)$, $f^2(\beta) = f(\alpha)  +2$,
then either for all $(\alpha, \beta) \in $ $P(\Sigma)$,
$(f(\alpha), f(\beta))$ \newline $= (-2, 0)$ or for all $(\alpha, \beta) \in$
$P(\Sigma)$, $(f(\alpha), f(\beta)) = (2, \pm 2)$.

(b) If there exists a pair $(\alpha, \beta) \in$ $P(\Sigma)$ so that
$f^2(\beta) \neq f(\alpha) +2$, then there exits a pair  $(\alpha',
\beta') \in $ $P(\Sigma)$ so that $f^2(\beta') \neq f(\alpha') +2$ and
one of $f^2(\alpha')$ or $f^2(\beta')$ is not 4. \rm

\it Proof. \rm  To prove (a), fix $\alpha \in  S'(\Sigma)$ which
bounds a $\Sigma_{1,1}$ and let  $P_{\alpha}$ be the
set of all non-separating classes  $\beta$ lying in $\Sigma_{1,1}$.
Take three elements $\beta_1, \beta_2$ and $\beta_3$ in $P_{\alpha}$
forming a triangle in the modular configuration. Let $f(\alpha) +2
$ be  $\mu^2$.  Then $f^2(\beta_i) = \mu^2$. By proposition 3.4(a),
we obtain $3\mu^2 \pm \mu^3 = \mu^2$. Thus either $\mu=0$ or 
$\mu^2= 4$, i.e., either $(f(\alpha), f(\beta)) = (-2, 0)$
for all $\beta \in P_{\alpha}$ or $(f(\alpha), f(\beta)) = (2, \pm 2)$
for all  $\beta \in P_{\alpha}$. 

To finish the proof of (a), we need to show that the above two
cases cannot occur simultaneously. The above proof shows
that being  $(f(\alpha), f(\beta)) = (2, \pm 2)$ or $(-2, 0)$
depends only on the 1-holed torus bounded by $\alpha$. Thus
part (a) follows from  corollary 3.8(b).

To prove part (b),  we may assume that the characteristic of the field
$K$ is not 2 (otherwise by part (a) the result follows). 
 Now suppose otherwise
that for all $(\alpha, \beta) \in $  $P(\Sigma)$ so that $f^2(\beta) \neq
f(\alpha) + 2$, we have $f^2(\alpha) = f^2(\beta) = 4$. This implies
that $f(\alpha) = -2$. Consider the level-1
subsurface $\Sigma_{1,1}$ bounded by $\alpha$ which contains $\beta$.
Let $\beta$, $\beta_2$ and $\beta_3$ be three classes in $S'(\Sigma_{1,1})$
which form a triangle  in the modular configuration. 
Then for $i=2,3$, 
either $f^2(\beta_i) = 4 $ (if  $f^2(\beta_i) \neq
f(\alpha) + 2$) or  $f^2(\beta_i) = 0$ (if 
 $f^2(\beta_i)  = f(\alpha) + 2$).
By proposition 3.4(a),
we have $f^2(\beta) + f^2(\beta_2) + f^2(\beta_3) - f(\beta) f(\beta_2) f(\beta_3) = f(\alpha) + 2$. Thus $4 + f^2(\beta_2) + f^2(\beta_3) = \pm 2f(\beta_2)
f(\beta_3)$. But this is impossible since either $f^2(\beta_i) = 4 $ or $0$
for $i=2,3$.
$\square$

6.4. Let $\partial \Sigma_{1,2}$ be $b_1$ and $b_2$. By proposition 6.2
and lemma 6.3, it remains to prove the  following.

{\bf Proposition.} \it Suppose $f$ is a $K$ valued trace function  on
$S(\Sigma_{1,2})$ so that either (a) for all $(\alpha, \beta) \in P(\Sigma_{1,2})$
$(f(\alpha), f(\beta)) = (-2, 0)$ or (b) for all
$(\alpha, \beta) \in P(\Sigma_{1,2})$  
$(f(\alpha), f(\beta) ) = (2, \pm 2)$. Then $f$ is a character. \rm

6.5. We construct a  representation whose character is $f$
satisfying condition (a) in the proposition 6.4 in this section.

{\bf Lemma.} \it Let  $\partial \Sigma_{1,2}$ be $ b_1 \cup b_2$. Under the assumption of  proposition 6.4(a), we have
$f^2(b_i) = 4$ and $f(b_1) f(b_2) = -4$.
\rm

\it Proof. \rm 
Take $(\alpha_1, \beta_1) \in P(\Sigma_{1,1})$ and let $\Sigma_{0,4}$
and $\Sigma_{1,1}$ be the subsurfaces bounded by $\beta_1$ and
$\alpha_1$. Take $\beta_2$ to be a non-separating class lying in
$\Sigma_{0,4}$ so that $\beta_2  \perp_0 \alpha_1$. Then
$\alpha_1 \beta_2$ and $\beta_2 \alpha_1$ are both non-separating.
By the assumption, $f(\beta_i) = f(\beta_2 \alpha_1) = f(\alpha_1
\beta_2) = 0$ and $f(\alpha) = -2$. By proposition 3.5(a) applied
to $\Sigma_{0,4}$ with respect to the triangle $(\alpha_1, \beta_2,
\alpha_1 \beta_2)$, we  obtain $f(b_1) + f(b_2) = 0$.  By
proposition 3.5(b) applied to $\Sigma_{0,4}$ with respect to
$\{ \alpha_1 = \beta_2(\alpha_1 \beta_2), (\alpha_1 \beta_2) \beta_2\}$, 
we obtain $f(b_1) f(b_2) = -4$. Thus the result follows.
$\square$

Here is a construction of a
 representation $\rho: \pi_1(\Sigma_{1,2}) \to SL(2,K)$
whose character is $f$.
 For simplicity, let $f(b_1) =2$ and $f(b_2) = -2$.  Let
$\Sigma_{1,1}$ be obtained by attaching a disc to the
$b_1$ boundary component of the surface $\Sigma_{1,2}$. By corollary 3.4,
there is a representation $\rho_0 : \pi_1(\Sigma_{1,1}) \to SL(2, K)$
so that $tr(\rho_0(\alpha)) = 0$ for all $\alpha \in S'(\Sigma_{1,1})$
and $tr(\rho_0(b_2)) = -2$. Define $\rho = \rho_0 \circ i$ where
$i: \pi_1(\Sigma_{1,2}) \to \pi_1(\Sigma_{1,1})$ is the homomorphism
induced by the inclusion map.  Since $i$ send non-separating classes
to non-separating classes, it follows that the character of $\rho$ is $f$.

6.6.
We construct a representation whose character is the 
trace function $f$ satisfying condition (b) in the proposition 6.4.

We may assume that the characteristic of the field $K$ is not 2 in
this section (otherwise it is covered by \S 6.5).  

{\bf Lemma.} \it 
Under the assumption 
of proposition 6.4 (b), we have $f(b_1) = f(b_2) = \pm 2$. In
particular, $f$ is reducible over all level-0 subsurfaces. \rm

\it Proof. \rm Since  $f(\alpha) = 2$ for all separating classes
$\alpha$, $f$ is reducible  on all 1-holed tori. In particular,
if $\alpha_1 \perp \alpha_2$, then by proposition 3.4(a), 
$f(\alpha_1 \alpha_2) = \frac{1}{2}f(\alpha_1) f(\alpha_2)$. Thus
$f(\alpha_1 \alpha_2) = f(\alpha_2 \alpha_1)$. On the other hand,
if $\beta$ and $\gamma$ are two non-separating classes so that 
$\beta \perp _0 \gamma$, then there exists three non-separating
classes $\delta_1 \perp \delta_2 \perp \delta_3$ so that
$\beta = \delta_1 \delta_2 \delta_3$ and $\gamma = \delta_3 \delta_2
\delta_1$. (Indeed, the pair ($\beta, \gamma)$ is unique up
to the homeomorphism of the surface). Thus we have $f(\beta)
= f(\gamma)$.

\midspace{0.1cm}
\centerline{\epsfbox{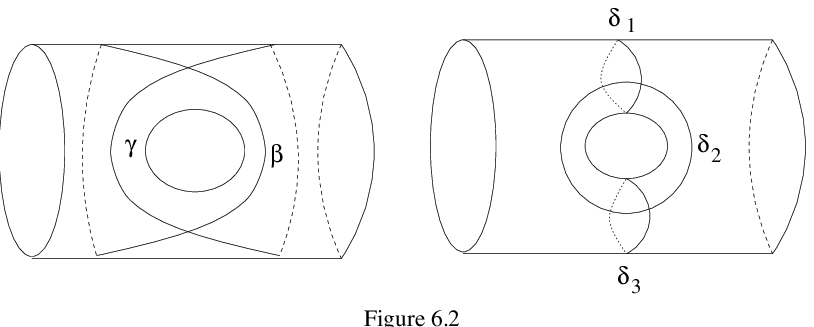}}

Now take two classes $\alpha_1 \perp \alpha_2$. Since $f(\alpha_1)
f(\alpha_2) f(\alpha_1 \alpha_2) = 8$, we may assume that
$f(\alpha_1) = 2$. Let $\Sigma_{0,4}$ be the subsurface bounded
by $\alpha_1$ and let $\alpha_3 \perp _0 \alpha_4$ be two
non-separating classes lying in $\Sigma_{0,4}$. Then both
$\alpha_5 = \alpha_3 \alpha_4$ and $\alpha_5' = \alpha_4 \alpha_3$
are separating classes. By the observation above, $f(\alpha_3)
= f(\alpha_4)$ and $f(\alpha_5) = f(\alpha_5') = 2$. By proposition
3.5(a) applied to the triangle $(\alpha_3, \alpha_4, \alpha_5)$ in
$\Sigma_{0,4}$,
we obtain $f(b_1) + f(b_2) = \pm 4$. By proposition 3.5(b)
applied to $\alpha_5, \alpha_5'$, we obtain $f(b_1) f(b_2) = 4$.
Thus $f(b_1) = f(b_2) = \pm 2$. Since $f(\alpha) =2$ for
all separating classes, this shows that $f$ is reducible on
all level-0 subsurface bounded by $\alpha$. Thus $f$ is
reducible on all level-0 subsurfaces.
$\square$

We now construct a diagonalizable representation $\rho$ of
$\pi_1(\Sigma_{1,2})$ whose character is $f$ as follows.

Take $(\alpha, \beta) \in P(\Sigma_{1,2})$ and let $X$ and $Y$
be the level-1 subsurfaces bounded by $\alpha$ and $\beta$ 
respectively. By the same argument as in \S6.2, we construct a
diagonalizable  $SL(2,K)$ representation of $\pi_1(\Sigma_{1,2})$
so that its character $g$ equals $f$ on
$S(X) \cup S(Y)$. To show that $g = f$, by corollary
3.7, it suffices to prove $f(\gamma) = g(\gamma)$ for
all $\gamma \perp_0 \alpha$ and $\gamma \perp \beta$.  Let
$\delta$ be a class disjoint from $\alpha$ and $\gamma$. Then
$f(\delta) = g(\delta)$.  By the reducibility of $f$ and $g$
on the level-0 subsurface bounded by $\delta$, $\gamma$,
it follows that $f(\gamma) = \frac{1}{2}f(b_1) f(\delta)
=\frac{1}{2}g(b_1) g(\delta) = g(\gamma)$.
$\square$

\S7. {\bf The Proof of Theorem 1.2}

The goal of this section is to prove theorem 1.2 for surfaces of
positive genus by using the induction on the level of the surface.

Let $\Sigma = \Sigma_{g,n}$ be a surface of positive genus
and $f$ a $K$-valued trace function defined on $S(\Sigma)$. Recall
that $P(\Sigma)$ is defined to be $\{(\alpha, \beta) \in S(\Sigma_{g,n})
\times S(\Sigma_{g,n}) |$ $\alpha$ bounds a $\Sigma_{1,1}$ and
$\beta$ is a non-separating class lying in the subsurface $\Sigma_{1,1}$\}.
The proof breaks into the following  two cases: (a)
there exists $(\alpha, \beta) \in $ $P(\Sigma)$ so that $f^2(\beta) \neq
f(\alpha) + 2$, and (b) for all $(\alpha, \beta) \in $ $P(\Sigma)$
$f^2(\beta) = f(\alpha) + 2$.  By lemma 6.3, case (b) is equivalent to
two subcases (b1) for all $(\alpha, \beta) \in $ $P(\Sigma)$, 
$(f(\alpha), f(\beta)) = (-2, 0)$ and (b2) for all
$(\alpha, \beta) \in $ $P(\Sigma)$, $(f(\alpha), f(\beta)) = (2, \pm 2)$.

We will deal with these three cases (a), (b1) and (b2) separately 
in the following sections.

7.1. Suppose the case (a) occurs. By lemma 6.3(b), we may assume that
one of $f^2(\alpha)$ or $f^2(\beta)$ is not 4. Let $X$ be the level-1 subsurface
bounded by $\alpha$ and let $Y$ be the subsurface $\Sigma -int(N(\beta))$.
Since $f|_{S(Y)}$ takes  some values other than $\pm 2$, by the
induction hypothesis if $g \geq 2$ and by the result in \S5 if
$g=1$, $f|_{S(Y)}$
is a character, say $f|_{S(Y)} = \chi_{\rho_Y}$ for an $SL(2, K)$
representation $\rho_Y$ of $\pi_1(Y)$.  Now if $f(\alpha) \neq 2$,  
then both $f|_{S(X)}$ and $f|_{S(X \cap Y)}$ are irreducible. Let
$\rho_X : \pi_1(X) \to SL(2,K)$ be any representation whose character
is $f|_{S(X)}$. By lemma 2.4, we may assume after a conjugation that
that $\rho_X |_{\pi_1(X \cap Y)} = \rho_Y|_{\pi_1(X \cap Y)}$. 
If $f(\alpha) = 2$, then both $f|_{S(X)}$ and $f|_{S(X \cap Y)}$
are reducible. Since one of $f^2(\alpha)$ or $f^2(\beta)$ is not
4, by lemma 2.6, there exist exactly two $SL(2,K)$ conjugacy
classes of representations of $\pi_1(X)$ whose characters
are $f|_{S(X)}$. Thus we may choose an $SL(2,K)$ representation
$\rho_X$ of $\pi_1(X)$ so that $\rho_X|_{\pi_1(X \cap Y)}
= \rho_Y|_{\pi_1(X \cap Y)}$ and $\chi_{\rho_X} = f|_{S(X)}.$
 
Define a representation $\rho: \pi_1(\Sigma) \to SL(2,K)$ by
$\rho|_{\pi_1(X)} = \rho_X$ and $\rho |_{\pi_1(Y)} = \rho_Y$.
Let $g$ be the character of $\rho$. Then
$g|_{S(X) \cup S(Y)}
= f|_{S(X) \cup S(Y)}$. 

To show that $g=f$, by corollary 3.7, it suffices to prove $f(\gamma)
=g(\gamma)$ for all $\gamma \perp_0 \alpha$ and $\gamma \perp \beta$.
Given such $\gamma$, consider the level-2 subsurface $\Sigma_{1,2}$
containing both $X$ and $\gamma$. Then by the proof of
 proposition 6.2 applied
to $\Sigma_{1,2}$ with respect to the decomposition $X$ and $Y \cap \Sigma_{1,2}$,
it follows that $g(\gamma) = f(\gamma)$.

7.2. To show the remaining cases, we need,

{\bf Lemma.} \it  Suppose $f$ is a trace function on $\Sigma_{g,r}$
 so that the case (b1) or (b2) holds. Then $f$ is reducible on
all level-0 subsurfaces. 

\rm

\it Proof. \rm 
Since each level-0 subsurface is contained in a 3-holed
torus subsurface, it suffices to prove the lemma for the
3-holed torus $\Sigma_{1,3}$.
Each level-0 subsurface in $\Sigma_{1,3}$ is either
 contained in a 2-holed
torus subsurface or is bounded by a boundary class. By \S6.5 and
\S6.6, those level-0 subsurfaces contained in a 2-holed torus
are reducible. It remains to show the reducibility of the level-0
subsurface $\Sigma_{0,3}(\gamma)$ bounded by a boundary class
$\gamma$ and two boundary components $b_1$ and $b_2$ of
$\Sigma_{1,3}$. Take disjoint non-separating classes $\gamma_1$
and $\gamma_2$ so that $\gamma_i \cap \gamma = \emptyset$ for $i=1,2$
and take $\gamma_3$ so that $\gamma_3 \cap \gamma_i = \emptyset$
for $i=1,2$ and $\gamma_3 \perp_0 \gamma$. Note that $\gamma_3$, $\gamma_3
\gamma$ and $\gamma \gamma_3$ are non-separating classes. 

In the  case (b1), $(f(\alpha), f(\beta)) = (-2, 0)$ for all $(\alpha,
\beta) \in P(\Sigma_{1,3})$. By proposition 3.5(a) applied
to the level-1 subsurface bounded by $\gamma_1$, $\gamma_2$ and
the triangle $(\gamma, \gamma_3, \gamma \gamma_3)$, we obtain
$f^2(b_1) + f^2(b_2) + f^2(\gamma) - f(\gamma) f(b_1) f(b_2) -4 =0$. Thus
by proposition 3.3, $\Sigma_{0,3}(\gamma)$ is reducible. 

In the case (b2), we may assume that the characteristic of $K$ is not
$2$ (otherwise the result is clear).  Now both $\Sigma_{0,3}(\gamma_3,
\gamma_i)$, $i=1,2$ and $\Sigma_{0,3}(\gamma_1, \gamma_2, \gamma)$
are reducible since they lie in some 2-holed torus subsurfaces. Thus,
by corollary 3.3, $f(b_i) = \frac{1}{2}f(\gamma_3) f(\gamma_i)$,
$i=1,2$ and $f(\gamma) = \frac{1}{2} f(\gamma_1)f(\gamma_2)$. This
implies that $f(b_1) f(b_2) f(\gamma) = 8$. By corollary 3.3, this
shows that $\Sigma_{0,3}(\gamma)$ is reducible.
$\square$

7.3. We now show that in the cases (b1) or (b2), the trace function
 $f$ is a character.

Take $(\alpha, \beta) \in $ $P(\Sigma)$
and let $X$ be the  $\Sigma_{1,1}$ subsurface bounded by $\alpha$ and
let $Y$ be the subsurface bounded by $\beta$. 

If $(f(\alpha), f(\beta)) = (2, \pm 2)$, then by lemma 7.2,
 we construct
a diagonalizable representation $\rho$ of $\pi_1(\Sigma)$ so that
its character $g = \chi_{\rho}$ satisfies $g|_{S(X) \cup S(Y)}
= f|_{S(X) \cup S(Y)}$. To show that $f=g$, by corollary 3.7,
it suffices to prove $f(\gamma) = g(\gamma)$ for all $\gamma \perp_0
\alpha$ and $\gamma \perp \beta$. Consider the level-2 subsurface
$\Sigma_{1,2}$ containing $X$ and $\gamma$.  Since both $f$ and $g$
are reducible on all level-0 subsurfaces,  by the proof of proposition 6.2,
it follows that $f(\gamma) = g(\gamma)$.

If $(f(\alpha), f(\beta)) = (2, 0)$ and the characteristic of
$K$ is not 2, we note that the genus $g$ of $\Sigma_{g,n}$ must
be 1. Indeed, if $g \geq 2$, then there exists essential subsurface
$\Sigma_{1,2}$ whose boundary components $\beta_i$, $i=1,2$
are non-separating simple
loops in $\Sigma_{g,n}$. By the assumption $f(\beta_i) = 0$. But by
lemma 6.5 applied to $\Sigma_{1,2}$, we have $f(\beta_i) = \pm 2$ which
is a contradiction.

We now construct a representation as follows. Let $\partial \Sigma_{1, n}$
be $b_1, ..., b_n$ and let $i : \pi_1( \Sigma_{1, n}) \to
\pi_1( \Sigma_{1,1})$ be the homomorphism induced by the inclusion map
$j: \Sigma_{1,n} \to \Sigma_{1,1}$ so that $j(b_n) = \partial \Sigma_{1,1}$.
Let $\rho': \Sigma_{1,1} \to SL(2,K)$ be a representation so that
$\chi_{\rho'}(\alpha) = 0$ for all non-separating class $\alpha$ and
$\chi_{\rho'}(\partial \Sigma_{1,1}) = -2$ (see corollary 3.4). 
Let $\rho_0 = \rho' \circ i$
be a representation of $\pi_1(\Sigma_{1,n})$. The fundamental
group  $\pi_1(\Sigma_{1,n})$
is a free group on $(n+1)$ generators $x_1, ..., x_{n+1}$ where
$x_1, ..., x_{n-1}$ correspond to the boundary components $b_1, ..., b_{n-1}$.
Now  modify $\rho_0$ to 
produce a new representation $\rho$ of $\pi_1(\Sigma_{1,n})$
by redefining $\rho(x_i) = \pm \rho_0(x_i)$  
so that $\chi_{\rho}(x_i) = f(b_i)$ for $i=1,2,..., n-1$.  Let $g$ be
the character of $\rho$ defined on $S(\Sigma_{1,n})$. Then
$g$ satisfies  $(g(\alpha), g(\beta) ) = (-2, 0)$ for all
$(\alpha, \beta) \in  P(\Sigma_{1, n})$ (indeed each non-separating
loop in $\Sigma_{1,n}$ becomes a non-separating loop in $\Sigma_{1,1}$).
Furthermore, by lemma 7.2, the character 
$g$ is reducible over all level-0 subsurfaces.
We prove that $f=g$ by induction on $n$. The result follows for $n=1,2$.
We first claim that $f(b_n) = g(b_n)$. To see this, take a boundary
class $\alpha'$ so that $\Sigma_{0,3}(\alpha')$ contains $b_n$ and $b_{n-1}$.
By the induction hypothesis applied to the subsurface $\Sigma_{1, n-1}$
bounded by $\alpha'$, we conclude that $f(\alpha') = g(\alpha')$.
By the reducibility of $f$ and $g$ on $\Sigma_{0,3}(\alpha')$ and
$f(b_{n-1}) = g(b_{n-1})$, it follows that $f(b_n) = g(b_n)$. Now for
any separating $\gamma$, let $\Sigma'$ be the planar subsurface
bounded by $\Sigma'$ in $\Sigma_{1,n}$.  Since $f$ and $g$ are
both reducible on all level-0 subsurfaces, $f$ and $g$ are reducible
on $\Sigma'$. Furthermore, $f$ and $g$ have the same values
on all but one boundary component $\gamma$ of $\Sigma'$. Thus,
by the reducibility, $f(\gamma) = g(\gamma)$.
$\square$ 

\S8. {\bf Proof of Theorem 1.1}

We begin with the following special case of theorem 1.1.

8.1. {\bf Proposition.} \it Suppose $K$ is a quadratically closed 
field and $f$ $:G \to K$ is a $K$-trace function defined on a finitely
generated group $G$. Then $f$ is the character of an $SL(2,K)$
representation of the group. \rm

\it Proof. \rm We first show that the result holds for $G=F_n$,
the free group on $n$ generators. Consider $F_n$ as the fundamental
group $\pi_1(\Sigma_{1,n-1})$ of the genus 1 surface with $n-1$
boundary components. Then by the work  of [Hel], $f$ induces a
trace function, denoted by $f'$, defined on $S(\Sigma_{1,n-1})$.
By theorem 1.2, there exists a representation $\rho$ of the
fundamental group $\pi_1(\Sigma_{1,n-1})$ to $ SL(2,K)$ whose
character is $f'$. Thus $\chi_{\rho}(x) = f(x)$ for each
$x \in \pi_1(\Sigma_{1, n-1})$ which has a simple loop representative.
Now by the remark following corollary 2.2, $f$ is the character of
$\rho$ on $G$.

 for $\pi_1(\Sigma_{1, n-1})$ so that all elements $x_{i_1} x_{i_2}... x_{i_k}$
 $1 \leq i_1 < i_2 < ..., < i_k \leq n$, $1 \leq k \leq n$, are represented by simple loops. \rm

 required property. Fig.

For the general n-generator group $G$, we follow an observation
of Gonz\`alez-Acu\~na and Montesinos-Amilibia [GM]. 
Let $\phi: F_n \to G$ be an epimorphism with $ker(\phi) = H$.
Then $g = f \circ \phi$ is a $K$-trace function defined on $F_n$.
Thus there exists a representation $\rho$ of $F_n$ whose character
is $g$. Furthermore, by the construction $tr \rho(x) = 2$ for all $x \in H$.
Now we use the following lemma of [GM].

{\bf Lemma} ([GM]). \it Suppose $\rho: F_n =<x_1, ..., x_n> \to SL(2,K)$ is
a representation and $x \in F_n$ so that $tr \rho(x ) =2$ and 
$tr( \rho([x, x_i]) =2$ for all $i$. Then either $\rho(x) =id$ or
$\rho$ is reducible. \rm

Indeed, if $\rho(x) \neq id$, then $\rho(x)$ has a unique eigenspace
in $K^2$. But $tr(\rho([x, x_i])=2$ shows that this eigenspace
is invariant under all $\rho(x_i)$. Thus $\rho$ is
reducible. 

By the lemma, if $\rho$ is irreducible, then $\rho(x) = id$ for all
$x \in H$. In particular, the representation $\rho$ induces
a representation $\rho'$ of $G$ to $SL(2,K)$ whose character is $f$.
If $\rho$ is reducible, we may replace $\rho$ by its diagonalization
$\rho'$ without changing the character. Now $tr(\rho'(x)) = 2$ if and
only if $\rho'(x) = id$. Thus the same argument goes through and we
construct a representation whose character is $f$. 
$\square$

8.2. We now prove theorem 1.1 for any group $G$. Let $f$ be a $K$-trace
function defined on $G$. We shall consider the following three cases:
(1) there exist $x, y \in G$ so that $f([x,y])  \neq 2$, (2)
for all $x, y \in G$, $f([x,y]) = 2$ but there exists $t \in G$ so that
$f(t) \neq \pm 2$, (3)   for all $x,y \in G$, $f([x,y]) =2 $ and
$f(x) = \pm 2$.

In the first case, consider the subgroup $<x,y>$ and restriction $f|_{<x,y>}$.
By lemma 2.3, there exists an irreducible representation $\rho_0$ of
$<x,y>$ whose character is $f|_{<x,y>}$. Given any element $z \in G$, 
consider the subgroup $<x,y,z>$ and the restriction $f|_{<x,y,z>}$.
By lemma 2.3, there exists a representation $\rho: <x,y,z> \to SL(2,K)$
so that its character is $f|_{<x,y,z>}$. Both $\rho|_{<x,y>}$
and $\rho_0$ have the same character and both are irreducible. By
lemma 2.4, we may assume after  conjugating  $\rho$ by an element in $SL(2,K)$ so
that $\rho|_{<x,y>} = \rho_0$. We denote this representation
by $\rho_z : <x,y,z> \to SL(2,K)$. Note that since $\rho_z$ is
irreducible, $\rho_z$ is unique. Now define a map $\mu: G \to
SL(2, K)$ by $\mu(z)= \rho_z(z)$. Clearly $tr(\mu(z)) = f(z)$ by
the construction. We claim that $\mu$ is a representation. Indeed,
given  $z_1, z_2  \in G$, 
consider the 4-generator subgroup $<x,y,z_1, z_2>$ and the
restriction $f|_{<x,y,z_1, z_2>}$. By proposition 8.1, there
exists a representation $\delta : < x,y,z_1, z_2> \to SL(2,K)$
whose character is $f|_{<x,y,z_1, z_2>}$.  By lemma 2.4, we 
may assume after conjugating by an element in $SL(2,K)$ that
$\delta|_{<x,y>} = \rho_0$. 
Thus we obtain $\rho_{z_i} = \delta|_{<x,y, z_i>}$ for $i=1,2$.
In particular this implies that $\mu(z_1 z_2) = \mu(z_1) \mu(z_2)$.

In the case (2), we consider the subgroup $<x, y> $ where $y=x$ so that
$f^2(x) \neq 4$. Let $\rho_0$ be a diagonal representation of $<x,y>$
whose character is $f|_{<x,y>}$. Note that the assumption $f([a,b]) =2$
implies the reducibility of the representations on all 2-generator
subgroup. We go through the same argument as in the previous
paragraph by taking all representations  $\rho_z,
\delta$ to be diagonalizable. Since $f^2(x) \neq 4$, by lemma 2.6, these
representations are unique. Thus the result follows.

Finally in the case (3), we have $f(x) = \pm 2$ and $f([x,y]) = 2$ for
all $x,y \in G$. If the characteristic of $K$ is $2$, then $f = 0$ and
$f$ is the character of the trivial representation. If the characteristic
of $K$ is not 2, then by lemma 2.2 (b), we obtain
 $f(xy) = f(x) f(y)/2$. Define a representation $\mu$ of
$G$ by $\mu(x) = \left(  \matrix f(x)/2 & 0 \\ 0 & f(x)/2  \endmatrix \right)
$. Then the character of $\mu$ is $f$.
$\square$.

\it Remark. \rm As the proof shows, theorem 1.1 follows as long
as one establishes theorem 1.1 for the free group on 4 generators.

\bigskip

\S9. {\bf Some Questions}

There are some questions arising  from the above considerations
concerning the finite presentations.
It is shown in [Lu1] and [Lu2] that
there exists a finite set
$F \subset S(\Sigma)$ so that  the Teichmuller space
(respectively the space of measured laminations) of $\Sigma$
is defined  by the restrictions of the length functions to $F$
subject to a finite set of polynomial equations supported in
level-1 subsurfaces.
The analogous question for the mapping class group of the
surface seems to be open. Namely, whether
the mapping class group
has a finite presentation whose generators are
finitely many  Dehn twists and whose relations (in these generators)
are supported
in level-1 subsurfaces.
 A recent work of Gervais [Ge] shows that one can find a
finite set of Dehn twists generating the mapping class group
so that the relations (in these generators) are supported
in level-3 subsurfaces.
Motivated by these, it is natural to ask if there exists
a finite set  $F \subset S(\Sigma)$ so that $SL(2,K)$ characters
are determined by their restrictions to $F$ subject to
polynomial equations supported in level-1 subsurfaces.
The proofs in \S4 and \S6 strongly suggest that the
answer is affirmative.  If the answer is positive,
it also implies that the character variety of any finitely
generated group can be defined by the restrictions of
the characters to a finite set of group elements subject
to polynomial equations supported in 3-generator subgroups.
The work of [GM] shows that the one can take the equations to be
supported in
5-generator subgroups.

There are several other related problems which seem to be
intersting.
The first question is that given a topological
group and a complex valued continuous trace function on the group, is
it the character of a continuous $SL(2, \bold C)$ representation of the
group? The second question is whether theorem 1.1 remains true for the
characters of $GL(n, \bold C)$ representations. To be more  precise,
suppose $f$ is a complex valued function defined on the
fundamental group of a surface so that the restriction of the
function to each level-1 subsurface group is a $GL(2, \bold C)$
character. Is there a $GL(n, \bold C)$ representation  $\rho$ of the
fundamental group so that $tr(\rho(x)) = f(x)$ for all $x$ lying
in some level-1 subsurface?
The third question is motivated
by Royden's theorem [Ro] for the Teichm\"uller spaces. Suppose
$\phi$ is an algebraic automorphism of the $SL(2,\bold C)$ character
variety of a surface group preserving the peripheral structure. Is $\phi$
induced by a self-homeomorphism of the surface?
Finally the analogous result to Jorgensen's discreteness criterion
seems to be the following. Suppose $\rho$ is a faithful representation
of a closed surface group to  $SL(2,\bold C)$ so that $\rho$
is discrete when restricted to each level-1 subsurface group. Is
$\rho$ discrete?

{\bf Appendix :  A Proof of Lemma 2.3}

{\bf Lemma 2.3}. \it Suppose $K$ is a field in which all quadratic equations with
coefficients in $K$ have roots in $K$. Given six numbers $x_1, x_2 , x_3, x_{12}, x_{23}$ and
$x_{31}$ in $K$, there exist three matrices $A_1, A_2,$ and $A_3$ in $SL(2, K)$ so that
$tr A_i = x_i$ and $trA_iA_j = x_{ij}$, for $i=1,2,3$ and $(i,j) = (1,2), (2,3),$
and $(3,1)$. \rm

\it Proof. \rm We divide the proof into three cases: in case 1, some $x_i \neq 
2$, in case 2, some $x_{ij} \neq \pm 2$, and in case 3, all $x_i$'s and 
$x_{ij}$'s are $\pm 2$.

\it Case 1. \rm Some $x_i \neq \pm 2$, say $x_1 \neq \pm 2$. Choose $\lambda $
in $K$ so that $x_1 = \lambda + \lambda^{-1}$. Clearly $\lambda \neq \pm 1$.
 Let $A_1 = \left( \matrix \lambda & 0 \\ 0 & \lambda^{-1}
\endmatrix \right)$, $A_2 = \left( \matrix a & b \\c & d \endmatrix \right)$, and
$A_3 = \left( \matrix x & y \\ z & w \endmatrix \right)$  be three $SL(2)$
matrices. We will find $a,b,c,d,x,y,z,w$ in $K$ solving the trace
equations. By $trA_2 = x_2$ and
$trA_1A_2 = x_{12}$, we obtain $a+d = x_2$ and $\lambda a + \lambda^{-1} d = x_{12}$.
Because $\lambda \neq \pm 1$, we can solve this system of linear equations uniquely
in $a, d$ in $K$. Similarly, by 
$trA_3 = x_3$ and $trA_3A_1 = x_{31}$,  we also solve $x, w$ uniquely in $K$.
It remains to find $b,c,y,z$ in $K$ so  that $bc = ad-1$, $yz= xw-1$ and $trA_2A_3 = x_{23}$,
i.e., $cy+bz = x_{23} -ax -dw$. If $ad -1 \neq 0$, i.e., $bc \neq 0$, choose
$b=1$. Let $c = ad-1$. Now, due to $bc \neq 0$, $cy + bz =  x_{23} -ax -dw$ and
$yz = xw -1$ can be solved in terms of $y,z$. 
If $ad-1 = 0$, there are two more subcases: $p = x_{23} - ax -dw \neq 0$ or $p=0$.
If $p=0$, we take $b=c=0$ and choose any pair $y, z$ so that $yz = xw-1$. If
$p \neq 0$, choose $b=1$, $c =0$. Then we have $z = p \neq 0$ and $y = (xw-1)/p$. 
Thus, in all cases, we find three matrices in $SL(2,K)$
satisfying the trace equations.

\it Case 2. \rm Some $x_{ij} \neq \pm 2$, say $x_{12} \neq \pm 2$. Then by case 1 
applied to the six ordered numbers $\{x_{12}, x_2, x_2x_3-x_{23}, x_1, x_{31}, 
x^2_2x_3- x_2x_{23} - x_3\}$,
we find three $SL(2,K)$ matrices $B_1, B_2, $ and $B_3$ 
so that $trB_1 = x_{12}$, $trB_2 = x_2$,
and $trB_3 = x_2x_3-x_{23}$, $trB_1B_2 = x_{1}$, $trB_1B_3 = x_{13}$
 and $tr B_2B_3 = x^2_2x_3  -x_2x_{23} -x_3$.
(Indeed, we take $B_1 = A_1A_2, B_2 = A_2^{-1}$ and $B_3 = A_2^{-1}A_3$ to find the six
numbers). Now let $A_1 = B_1B_2, A_2 = B_2^{-1}$ and $A_3 = B_2^{-1} B_3$. By the
basic trace identity (lemma 2.2(a)), it follows that $trA_i = x_i$ and $trA_iA_j = x_{ij}$.

\it Case 3. \rm All $x_i$'s and $x_{ij}$'s are $\pm 2$. First we note that if
$trA_i = x_i$ and $trA_iA_j = x_{ij}$, then $(-A_1, A_2, A_3)$ solves the
problem for the six numbers $\{-x_1, x_2, x_3, -x_{12}, x_{23}, -x_{31}\}$. Thus,
by changing the signs of $x_i$'s if necessary, we may assume that $x_1 = x_2 = x_3 = 2$.
There are four cases for $(x_1, x_2, x_2, x_{12}, x_{23}, x_{31})$ up to
symmetry: $(2,2,2,2,2,2)$, $(2,2,2,-2, 2,2)$, $(2,2,2,-2,-2, 2)$ and
$(2,2,2,-2, -2, -2)$. The corresponding solutions are listed below.

For $(2,2,2,2,2,2)$,  a solution $(A_1, A_2, A_3) $ is $ (id, id, id)$. For
$(2,2,2,-2, 2, 2)$, a solution is $( \left( \matrix 1 & 1 \\0 & 1 
\endmatrix \right),$$ \left( \matrix 1 & 0 \\-4 & 1 \endmatrix \right),
$$ \left( \matrix
1 &0 \\0 & 1 \endmatrix \right)) $. 
For $(2,2,2,-2,-2,2)$, a solution is 
$( \left( \matrix 1& 1 \\0 & 1 \endmatrix \right),$
$ \left( \matrix 1 & 0 \\-4 & 1 \endmatrix \right),$
$ \left( \matrix 1 & 0 \\ -4 & 1 \endmatrix \right))$. 
Finally for $(2,2,2,-2,-2,-2)$, a solution is
$(\left( \matrix 1 & 1 \\0 & 1 \endmatrix \right),$
$\left(\matrix 1 & 0 \\ -4 & 1 \endmatrix \right),$
$\left(\matrix -1 & 1 \\-4 & 3 \endmatrix \right))$.
$\square$

\centerline{ \bf References}

[Bu] Buser, P.: Geometry and spectra of compact Riemann surfaces. Progress
in Mathematics. Birkh\"auser, Boston, 1992.

[BG] Bers, L.; Gardiner, Frederick P.:
 Fricke spaces. Adv. in Math. 62 (1986), no. 3,
249-284. 

[BH] Brumfiel, G. W.; Hilden, H. M.: $SL(2)$ representations of finitely presented groups.
Contemporary Mathematics, 187. American Mathematical Society, Providence, RI, 1995. 

[CS] Culler, M.; Shalen, P.: Varieties of group representations and
splittings of 3-manifolds. Ann. Math. 117, (1983), 109-146.

[De] Dehn, M.: Papers on group theory and topology. J. Stillwell (eds.).
 Springer-Verlag, Berlin-New York, 1987.

[FK] Fricke, R.; Klein, F.: Vorlesungen  \"uber die Theorie der
Automorphen Functionen. Teubner, Leipizig, 1897-1912.

[Ge] Gervais, S.: A finite presentation of the
mapping class group of an orientable surface. preprint, 1998.

[Go] Goldman, W.: Topological components of spaces of representations.
Invent. Math. 93, (3) (1988), 557-607.

[Gr] Grothendieck, A.:  Esquisse d'un programme.
London Math.  Soc. Lecture Note Ser., 242, Geometric Galois actions, 
1, 5-48, Cambridge Univ.  Press, Cambridge, 1997.

[GM]
 Gonz\`alez-Acu\~na, F.; Montesinos-Amilibia, J.:
 On the character variety of group
representations in $SL(2, \bold C )$ and $PSL(2,\bold C)$. Math. Z. 214 (1993)
, no. 4, 627-652.

[Har] Harer, J.:
Stability of the homology of the mapping class groups of
orientable surfaces. Ann. of Math. (2) 121 (1985), no. 2, 215-249.

[Hav] 
Harvey, W. J.: Boundary structure of the modular group. Riemann surfaces and related
topics: Proceedings of the 1978 Stony Brook Conference (State Univ. New York, Stony Brook, N.Y., 1978), pp. 245-251,
Ann. of Math. Stud., 97, Princeton Univ. Press, Princeton, N.J., 1981.

[Hel] Helling, H.: Diskrete Untergruppen von SL(2, $\bold R$).
 Inventiones Math.,
17, (1972), 217-229.

[Hem] 
Hempel, J.:  Residual finiteness for $3$-manifolds. Combinatorial group theory and
topology (Alta, Utah, 1984), 379--396, Ann. of Math. Stud., 111, Princeton Univ. Press, Princeton, NJ, 1987.

[Ho] Horowitz, R.: Characters of free groups represented in the 2-dimensional
special linear group. Comm. Pure and App. Math.  25, (1972), 635-649.

[Jo] 
Jorgensen, T.: On discrete groups of M\"obius transformations. Amer. J. Math. 98 (1976),
no. 3, 739-749.

[Ke] Keen, L.: On Fricke moduli. 
1971 Advances in the Theory of Riemann Surfaces (Proc. Conf.,
Stony Brook, N.Y., 1969) pp. 205--224 Ann. of Math. Studies, No. 66. Princeton Univ. Press, Princeton, N.J. 

[Li1] Lickorish, R.: A representation of oriented combinatorial 3-manifolds. Ann.
Math. 72,  (1962), 531-540.

[Li2] Lickorish, R.: private communication.

[Lu1] Luo, F.:  Geodesic length functions and Teichm\"uller spaces.
J. Differential Geom.   48, (1998), 275-317.

[Lu2] Luo, F.: Simple loops on surfaces and their intersection numbers,
preprint, 1997. The Research announcement appeared in  Math. Res. Letters, 
 5, (1998), 47-56.

[Lu3] Luo, F.:  Automorphisms of the curve complex. Topology, 
to appear.

[Lu4] Luo, F.:  Grothendieck's Reconstruction Principle and 2-Dimensional 
Topology and Geometry. Communications in Contemporary Mathematics, to appear.

[LM]
Lubotzky, A.; Magid, A.: Varieties of representations of finitely generated
groups. Mem. Amer. Math. Soc. 58 (1985), no. 336.

[Ma] Magnus, W.: Rings of Fricke characters and automorphism groups of
free groups. Math. Zeit. 170, (1980), 91-103.

[MS] Moore, C.; Seiberg, N.: Polynomial equations for rational
conformal field theories. Phys. Lett. B.  212, (1988), 451-460.

[Pr]
Procesi, C.: Finite dimensional representations of algebras. Israel J. Math. 19 (1974),
169-182.

[Ro]
 Royden, H. L.: Automorphisms and isometries of Teichm\"uller space. 1971 Advances in the
Theory of Riemann Surfaces (Proc. Conf., Stony Brook, N.Y., 1969) pp. 369-383 Ann. of Math. Studies, No. 66. Princeton Univ. Press, Princeton, NJ.

[Sa]
Saito, K.: Character variety of representations of a finitely generated group in $ SL_2$. Topology and Teichm\"uller spaces (Katinkulta, 1995), 253-264, World Sci. Publishing, River Edge, NJ, 1996.

[Th] Thurston, W.: On the geometry and dynamics of diffeomorphisms of
surfaces. Bul. Amer. Math. Soc.  19, (1988) no 2, 417-438.       

[Vo] Vogt, H.: Sur les invariants fondamentaux des \'equations 
diff\'erentielles lin\'eaires du second ordre. Ann. de l'Ecole Normale Superieur (3), 6, Suppl.
(1889), 3-72.

Department of Mathematics, Rutgers University, New Brunswick, NJ 08903

email: fluo\@math.rutgers.edu

\end

\end